\documentclass[12pt,reqno]{amsart}

\usepackage[arrow,matrix,curve]{xy}

\usepackage[dvips]{graphicx} 

\usepackage{amssymb, latexsym, amsmath, amscd, array,
}

\newtheorem{theorem}{Theorem}[section]

\newtheorem{principle}[theorem]{Principle}

\theoremstyle{definition}

\numberwithin{equation}{section}
\newcommand\R {{\mathbb R}}

\newcommand\parisotes{{$\pi\alpha\rho\iota\sigma\acute{o}\tau\eta\varsigma$}}

\newcommand\parisos{{$\pi\acute\alpha\rho\iota\sigma{o}\varsigma$}}

\newcommand\Stromholm{{Str\o mholm}}

\newcommand\Los{{\L}o{\'s}}

\title[Adequality from Diophantus to Fermat]{Almost Equal: The Method
of Adequality from Diophantus to Fermat and Beyond}

\author{Mikhail G. Katz}

\address{M.~Katz, Department of Mathematics, Bar Ilan University,
Ramat Gan 52900 Israel} \email{katzmik@macs.biu.ac.il}

\author{David M. Schaps}

\address{D. Schaps, Department of Classical Studies, Bar Ilan
University, Ramat Gan 52900 Israel}

\email{dschaps@mail.biu.ac.il}

\author{Steven Shnider} 

\address{S.~Shnider, Department of Mathematics, Bar Ilan University,
Ramat Gan 52900 Israel} \email{shnider@macs.biu.ac.il}

\begin{document}

\begin{abstract}
We analyze some of the main approaches in the literature to the method
of `adequality' with which Fermat approached the problems of the
calculus, as well as its source in the \parisotes{} of Diophantus, and
propose a novel reading thereof.

Adequality is a crucial step in Fermat's method of finding maxima,
minima, tangents, and solving other problems that a modern
mathematician would solve using infinitesimal calculus.  The method is
presented in a series of short articles in Fermat's collected works
\cite[p.~133-172]{Tan}.  We show that at least some of the
manifestations of adequality amount to variational techniques
exploiting a small, or infinitesimal, variation~$e$.  

Fermat's treatment of geometric and physical applications suggests
that an aspect of approximation is inherent in adequality, as well as
an aspect of smallness on the part of~$e$.  We question the relevance
to understanding Fermat of 19th century dictionary definitions of
\parisotes{} and \emph{adaequare}, cited by Breger, and take issue
with his interpretation of adequality, including his novel reading of
Diophantus, and his hypothesis concerning alleged tampering with
Fermat's texts by Carcavy.  We argue that Fermat relied on Bachet's
reading of Diophantus.  

Diophantus coined the term \parisotes{} for mathematical purposes and
used it to refer to the way in which 1321/711 is approximately equal
to 11/6.  Bachet performed a semantic calque in passing from
\emph{pariso\=o} to \emph{adaequo}.

We note the similar role of, respectively, adequality and the
Transcendental Law of Homogeneity in the work of, respectively, Fermat
\cite{Fer} and Leibniz \cite{Le84} on the problem of maxima and
minima.
\end{abstract}


\maketitle

\tableofcontents

\section{The debate over adequality}
\label{one}

Adequality, or~\parisotes{} (parisot\=es) in the original Greek of
Diophantus,%
\footnote{{\em adaequalitas\/} or {\em adaequare} in Latin.  See
Section~\ref{ety} for a more detailed etymological discussion.}
is a crucial step in Fermat's method of finding maxima, minima,
tangents, and solving other problems that a modern mathematician would
solve using infinitesimal calculus.  The method is presented in a
series of short articles in Fermat's collected works
\cite[pp.~133-172]{Tan}.  The first article, {\em Methodus ad
Disquirendam Maximam et Minimam\/},%
\footnote{In French translation, {\em M\'{e}thode pour la Recherche du
Maximum et du Minimum}.  Further quotations will be taken from the
French text \cite{Fer}.}
opens with a summary of an algorithm for finding the maximum or
minimum value of an algebraic expression in a variable~$A$.  For
convenience, we will write such an expression in modern functional
notation as~$f(a)$.%
\footnote{Following the conventions established by Vi\`ete, Fermat
uses capital letters of vowel, A, E, I, O, U for variables, and
capital letters of consonants for constants.}

\subsection{Summary of Fermat's algorithm}
\label{summary}

The algorithm can be broken up into six steps in the following way:

\begin{enumerate}
\item
Introduce an auxiliary symbol~$e$, and form~$f(a+e)$;
\item
Set {\em adequal\/} the two expressions~$f(a+e) =_{\text{AD}} f(a)$;%
\footnote{The notation ``$=_{\text{AD}}$'' for adequality is ours, not
Fermat's.  There is a variety of differing notations in the
literature; see Barner \cite{Bar} for a summary.}
\item
\label{cancel}
Cancel the common terms on the two sides of the adequality.  The
remaining terms all contain a factor of~$e$;
\item
\label{divide}
Divide by~$e$ (see also next step);
\item
\label{higher}
In a parenthetical comment, Fermat adds: ``or by the highest common
factor of~$e$";
\item
\label{among}
Among the remaining terms, suppress%
\footnote{\label{discard1}The use of the term ``suppress'' in
reference to the remaining terms, rather than ``setting them equal to
zero'', is crucial; see footnote~\ref{discard2}.}
all terms which still contain a factor of~$e$.%
\footnote{\label{meaningless1}Note that division by~$e$ in
Step~\ref{divide} necessarily \emph{precedes} the suppression of
remaining terms that contain~$e$ in Step~\ref{among}.  Suppressing all
terms containing~$e$, or setting~$e$ equal to zero, at
stage~\ref{cancel} would be meaningless; see
footnote~\ref{meaningless2}.}
Solving the resulting equation for~$a$ yields the extremum of~$f$.
\end{enumerate}

In modern mathematical language, the algorithm entails expanding the
difference quotient
\[
\frac{f(a+e)-f(a)}{e}
\]
in powers of~$e$ and taking the constant term.%
\footnote{Fermat also envisions a division by a higher power of~$e$ as
in step~\eqref{higher} (see Section~\ref{three}).}
The method (leaving aside step~\eqref{higher} for the moment) is
immediately understandable to a modern reader as the elementary
calculus exercise of finding the extremum by solving the
equation~$f'(a)=0$.  But the real question is how Fermat understood
this algorithm in his own terms, in the mathematical language of his
time, prior to the invention of calculus by Barrow, Leibniz, Newton,
et al.

There are two crucial points in trying to understand Fermat's
reasoning: first, the meaning of ``adequality'' in step (2), and
second, the justification for suppressing the terms involving positive
powers of~$e$ in step~\eqref{among}.  The two issues are closely
related because interpretation of adequality depends on the conditions
on~$e$.  One condition which Fermat always assumes is that~$e$ is
positive. He did not use negative numbers in his calculations.

Fermat introduces the term \emph{adequality} in \emph{Methodus} with a
reference to Diophantus of Alexandria.  In the third article of the
series, {\em Ad Eamdem Methodum\/} ({\em Sur la M\^{e}me
M\'ethode\/}), he quotes Diophantus' Greek term~\parisotes, which he
translates following Xylander and Bachet, as {\em adaequatio\/} or
{\em adaequalitas\/} \cite[p.~126]{Fer} (see A.~Weil
\cite[p.~28]{We84}).

\subsection{Etymology of \parisotes}
\label{ety}

The Greek word \parisotes{} consists of the prepositional prefix
\emph{para} and the root \emph{isot\=es}, ``equality".  The prefix
\emph{para}, like all ``regular" prefixes in Greek, also functions as
a preposition indicating position; its basic meaning is that of
proximity, but depending upon the construction in which it appears, it
can indicate location (``beside"), direction (``to"), or source
(``from") (see Luraghi \cite[p.~20-22, 131-145]{Lu}).

Compounds with all three meanings are found.  Most familiar to
mathematicians will be \emph{parall\=elos}
($\pi\alpha\rho\acute\alpha\lambda\lambda\eta\lambda o\varsigma$),
used of lines that are ``next to" one another; the Greek for ``nearly
resembling" is \emph{parapl\=esios}
($\pi\alpha\rho\alpha\pi\lambda\acute\eta\sigma\iota o\varsigma)$; but
we also find direction in words like \emph{paradosis}
($\pi\alpha\rho\acute\alpha\delta{}o\sigma\iota\varsigma$),
``transmitting, handing over" (from \emph{para} and \emph{dosis},
``giving"), and source in words like \emph{paral\=epsis}
%
%paralepsis, paraleksis
%
($\pi\alpha\rho\acute\alpha\lambda\eta\psi\iota\varsigma$), ``receiving"
(from \emph{para} and \emph{l\=epsis}, ``taking"); there are other
meanings for this prefix that do not concern us.

The combination of \emph{para} and \emph{isos} (``equal") can refer
either to a simple equality (Aristotle, \emph{Rhetoric} 1410a 23:
$\pi\alpha\rho\acute\iota\sigma\omega\sigma\iota\varsigma$~$\delta$'
\, '\!\!$\varepsilon\grave\alpha\nu$ \, '\!\!$\acute\iota\sigma\alpha$
$\tau\grave\alpha$~$\kappa \tilde\omega \lambda\alpha$, ``It is
\emph{paris\=osis} if the members of the sentence are equal") or an
approximate equality (Strabo 11.7.1, \,
`\!\!\!$\acute{\omega\;}\!\varsigma$
%
%os
%
$\varphi\eta\sigma\iota$~$\Pi\alpha\tau\rho
o\kappa\lambda\tilde\eta\varsigma$,\, 
%
%patrokles
%
`\!$\grave o \varsigma\,$%
%
%another os
%
%diacritic over ``o'' is unclear
%
$\kappa\alpha\grave\iota$
%
%kai
%
$\pi\acute\alpha\rho\iota\sigma o\nu$
%
%parison
%
\, `\!\!\!$\eta\gamma\varepsilon\tilde\iota\tau\alpha\iota$
%
%egeitai
%
$\tau \grave o$
%
%to
%
$\pi\acute\varepsilon\lambda\alpha\gamma o\varsigma$
%
%pelagos
%
$\tau o\tilde\upsilon\tau o$
%
%touto
%
$\tau\tilde{\displaystyle\mathop{\omega}_{'}}$
%
%to
%
$\Pi o\nu\tau\iota\kappa\tilde{\displaystyle\mathop{\omega}_{'}}$
%
%Pontiko
%
``As Patrocles says, who also considers this sea [the Caspian] to be
approximately equal [\emph{parison}] to the Black Sea"). We know of no
passage other than those of Diophantus in which a term involving
\emph{para} and \emph{isos} refers to mathematical equality, whether
approximate or otherwise.  Diophantus himself used the term
\emph{parisos} to describe terms that are approximately equal, as we
shall demonstrate below (Subsection~\ref{61}).

The term \emph{isot\=es} denotes a relationship (``equality"), not an
action (``setting equal"); the normal term for the action of
equalizing would be a form in \emph{-\=osis}, and in fact the words
\emph{is\=osis}
%
%isosis
%
('\!$\acute\iota\sigma\omega\sigma\iota\varsigma$) and
\emph{paris\=osis}
($\pi\alpha\rho\acute\iota\sigma\omega\sigma\iota\varsigma$)
%
%parisosis
%
are attested with the meaning ``making equal" or ``equalization"; the
latter is a common term in rhetoric for using the same (or nearly the
same) number of words or syllables in parallel clauses. The word
$\pi\alpha\rho\iota\sigma \acute o\tau\eta\varsigma$, which occurs
only in the two Diophantus passages, is on the face of it more
appropriate to the meaning ``near equality".

Fermat himself may not have gotten this far into Greek etymology. On
the other hand, Fermat viewed Diophantus through the lens of Bachet's
analysis.  Bachet does interpret it as approximate equality.  If
Fermat follows Bachet, seeking to interpret Fermat's method of
adequality/\parisotes{} based on the Latin term, \emph{adaequare}, is
missing the point (see Subsection~\ref{philo}).

\subsection{Modern interpretations}

There are differing interpretations of Fermat's method in the
literature.  A.~Weil notes that Diophantus uses the Greek term
\begin{quote}
to designate his way of approximating a given number by a rational
solution to a given problem (cf.~e.g.~{\em Dioph.\/}V.11 and 14) (Weil
1984 \cite[p.~28]{We84}).
\end{quote}
According to Weil's interpretation, approximation is implicit in the
meaning of the original Greek term.

H.~Breger rejects Weil's interpretation of Diophantus, and proposes
his own interpretation of the mathematics of Diophantus' \parisotes.
He argues that \parisotes{} means {\em equality\/} to Diophantus
(Breger \cite[p.~201]{Bre94}).

Thus, the question of whether there is an element of approximation in
the Greek source of the term is itself subject to dispute.  There is
also a purely algebraic aspect to adequality, based on the ideas of
Pappus of Alexandria, and Fermat's predecessor Vi\`ete.  In {\em Sur
la M\^{e}me M\'ethode} following the comment quoted above, Fermat
writes as follows:
\begin{quote} 
En cet endroit, Pappus appelle un rapport minimum
$\mu o \nu\alpha\chi \grave{o}\nu$~$\kappa\alpha \grave{\iota}$\,
'\!\!$\epsilon\lambda\acute{\alpha}\chi\iota\sigma\tau o \nu$
(singulier et minimum), parce que, si l'on propose une question sur
les grandeurs donn\'ees, et qu'elle soit en g\'en\'eral satisfaite par
deux points, pour les valeurs maximum et minimum, il n'y aura qu'un
point qui satisfasse.  C'est pour cela que Pappus appelle {\em minimum
et singulier\/} (c'est-\`a-dire unique) le plus petit rapport de tous
ceux qui peuvent \^etre propos\'es dans la question (Fermat
\cite[p.~127]{Fer}).
\end{quote}

The point is that the extremum of a quadratic expression at the
point~$a$ corresponds to a double root in~$e$ for what would be in
modern terms the equation~$f(a+e)-f(a)=0$.  From this point of view,
Fermat explains his method in terms of roots of algebraic equations.
In the first paragraph of the fourth article in the series, {\em
Methodus de Maxima et Minima} ({\em M\'ethode du Maximum et Minimum}),
Fermat reveals the source of this procedure:
\begin{quote}
En \'etudiant la m\'ethode de la {\em syncrise} et de la {\em
anastrophe} de Vi\`ete, et en poursuivant soigneusement son
application \`a la recherche de la constitution des \'equations
corr\'elatives, il m'est venu \`a l'esprit d'en d\'eriver un
proc\'ed\'e pour trouver le maximum et le minimum et pour r\'esoudre
ainsi ais\'ement toutes les difficult\'es relatives aux conditions
limites, qui ont caus\'e tant d'embarras aux g\'eom\`etres anciens et
modernes (Fermat \cite[p.~131]{Fer}).
\end{quote}

From this point of view, adequality is based on replacing the
variable~$a$ by the variable~$a+e$ in the original algebraic
expression and thus creating an equation in~$a$ and~$e$ which is
required to have a double root at~$e=0$ for an extremal point~$a$.
This interpretation is considered by Breger \cite{Bre94} and
K.~Barner~\cite{Bar} to cover all the examples.  They deny that any
kind of ``approximation" is involved, and hold adequality to be a
formal or algebraic procedure of ``setting equal".

\subsection{Wieleitner, \Stromholm, and Giusti}
\label{13}

The authors H.~Wieleitner \cite{Wi}, P.~\Stromholm~\cite{Strom}, and
E.~Giusti \cite{Giu} argue that both interpretations, algebraic and
approximation, are valid, representing different stages in the
development of Fermat's method.  The algebraic approach, following
Pappus and Vi\`ete, involves equating two values~$f(a)$ and~$f(a+e)$,
below the maximum or above the minimum.%
\footnote{\label{sometimes}Sometimes the second point is denoted~$e$
rather than~$a+e$, and the two expressions that are equated are~$f(a)$
and~$f(e)$.}
However, there is another point of view in which~$f(a)$ and~$f(a+e)$
are definitely not assumed by Fermat to be equal, as he writes in {\em
Sur la M\^{e}me M\'ethode} that he compares the two expressions
``comme s'ils \'etaient \'egaux, quoiqu'en fait ils ne le soient
point" (Fermat \cite[p.~126]{Fer}), and a little later ``une
comparaison feinte ou une {\em ad\'egalit\'e\/}" (ibid., p.~127]).
Giusti \cite{Giu} and \Stromholm~\cite{Strom} consider this to
represent a second stage in the development of Fermat's method.
%Revised: added the last sentence
Similarly, in the fourth article, which begins with a reference to
Vi\`ete and emphasizes the algebraic approach, Fermat introduces an
element of approximation, remarking that the difference between the
two points~$a$ and~$e$ goes to zero:%
\footnote{\label{kinetic}See footnote~\ref{sometimes}.}
\begin{quote}
Plus le produit des segments augmentera, plus au contraire diminuera
la diff\'erence entre~$a$ et~$e$, jusqu'\`a ce qu'elle s'\'evanouisse
tout \`a fait \cite[p.~132]{Fer}.
\end{quote}
What complicates the task is that Fermat does not separate clearly
between these two methods.  They appear in successive paragraphs.

M.~Mahoney understands one of the meanings of adequality as
``approximate equality'' or ``equality in the limiting case'' (Mahoney
\cite[p.~164, end of footnote~46]{Mah73}), while emphasizing that the
term has multiple meanings.

Fermat never gave a full explanation of his method, but he derived it
from three sources: Diophantus, Pappus, and Fermat's predecessor
Vi\`ete (Vieta).  If we consider the source in Diophantus and
interpret adequality as approximate equality as did Weil, \Stromholm,
and Giusti, it is natural to ask whether Fermat considered~$e$ to be
arbitrarily small and eventually negligible, although he never
explicitly stated such an assumption.  On the other hand, the
algebraic point of view, finding a condition for a unique root of
multiplicity 2, following Pappus and Vi\`ete, is clearly the point of
view in a number of examples mentioned above.

These considerations do not resolve the issue of what Fermat thought
about the actual magnitude of~$e$.

\Stromholm~\cite{Strom} and Wieleitner \cite{Wi} deal with this
question.  They distinguish two methods in Fermat.  One method is
algebraic, following Pappus and Vi\`ete, which \Stromholm~calls M2.
The other method, M1, is interpreted as expanding~$f(a+e)-f(a)$ in
powers of~$e$.  The latter approach is most fully expounded in
Fermat's letter to Br\^ulart~\cite{Fer2}.  The letter attempts to
explain why the method guarantees a maximum or minimum without
assuming a condition on the size of~$e$.

The ``approximation" interpretation actually branches out into two
distinct approaches.  The difference between them concerns the
interpretation of the symbol~$e$ that Fermat uses, to form expressions
that in modern notation would be written as ``$f(x+e)-f(x)$''.
Namely, one can think of~$e$ as representing a kinetic process such as
``tending to 0", as in the fourth article, as cited above, or one can
think of~$e$ as ``infinitesimal".  G.~Cifoletti \cite{Ci} and
J.~Stillwell \cite{Sti} interpret it in accordance with the latter
approach.

\subsection{Three approaches to the nature of~$E$}

There are therefore at least three different approaches to Fermat's
symbol~$E$ as it appears in adequalities: algorithmic or
formal/algebraic; kinetic~$E\to 0$; and infinitesimal.  We will argue
that the last one is closest to Fermat's thinking.

Finally, we note that two distinct issues are sometimes conflated in
the literature on Fermat's method.  The first issue is whether
adequality means ($\alpha$) ``setting equal", or whether it is
($\beta$) an ``approximate equality''.  A second, separate issue
concerns the question of what the famous symbol~$E$ stands for: is it
(A) an arbitrary variable, or does it imply (B) some notion of size:
small, infinitesimal, tending to zero, etc.

The trend in the literature is that scholars following the
interpretation ($\alpha$), also adopt (A), and similarly for the other
pair.  For instance, Breger \cite[p.~206]{Bre94} rejects the
small/infinitesimal idea and supports~(A).  But the thrust of his {\em
argument\/} is to support the~$(\alpha)$-interpretation rather than
the (A)-interpretation.

These are, in fact, separate issues, as can be seen most readily in
the context of an infinitesimal-enriched ring such as the dual numbers
$\mathcal D$ of the form~$a+b\epsilon$ where~$\epsilon^2=0$.  To
differentiate a polynomial~$p(x)$, we apply the following purely
algebraic procedure: expand in powers of~$\epsilon$; write
$p(x+\epsilon)=p(x) +\epsilon q(x)$, where~$q(x)$ has no~$\epsilon$
terms, then one has~$p(x+\epsilon)-p(x)=\epsilon q(x)$ and~$q(x)$ is
the derivative.

This produces the derivative of~$p$ over~$\mathcal D$.  A similar
procedure works over any other reasonable infinitesimal-enriched
extension of~$\R$ such as the hyperreals.  The algebraic nature of
this procedure in no way contradicts the infinitesimal nature
of~$\epsilon$.

In addition to arguing that ``adequality" is a purely algebraic
procedure of ``setting equal", Breger {\em believes\/} that~$E$ is a
formal variable with no assumption on size.  However, one does not
necessarily imply the other.  The procedures in non-standard analysis
are purely algebraic and can be programmed by a finite algorithm (no
need for an infinite limiting process), yet here~$E$ is definitely
infinitesimal.

No analysis of Fermat's method can be considered complete that does
not include a discussion of the application to transcendental curves.
Such an analysis appears in Section~\ref{cycloid}.

Concerning the question as to whether Fermat's method is a purely
algorithmic/algebraic one, with~$E$ being a formal variable, or
whether it involves some notion of ``smallness" on the part of~$E$, we
argue that the answer depends on which stage of Fermat's method one is
dealing with.  He certainly did present an algorithmic outline of his
method in a way that suggests that~$E$ is a formal variable.  However,
when one examines other applications of the method, one notices
additional aspects of Fermat's method which cannot be accounted for by
means of a ``formal" story.  Thus, Fermat exploits his adequality to
solve a least time variational problem for the refraction of light
(see Section~\ref{nine}).  Here~$E$ corresponds to a variation of a
physical quantity, and it would seem paradoxical to describe it as a
formal variable in this context.  Furthermore, in the case of the
cycloid, the transcendental nature of the problem creates a situation
that cannot be treated algebraically at all (see
Section~\ref{cycloid}).

Once it is accepted that at least in some applications, the aspect of
``smallness" on the part of~$E$ is indispensable, one can ask in what
sense precisely is~$E$ ``small".  Today we know of two main approaches
to ``smallness", namely, (1) by means of kinetic ideas related to
limits, or (2) by means of infinitesimals.  The former ideas were as
yet undeveloped in Fermat's time (though they are already present in
Newton only a few decades later), and in fact one finds very little
``tends to\ldots" material in Fermat.  Meanwhile, infinitesimals were
already widely used by Kepler, Wallis, and others.  Fermat's 1657
letter to Digby on Wallis's method (see Section~\ref{method}) shows
that he was intimately familiar with the method of indivisibles.  What
we argue therefore is that it is more reasonable to assume (2).

The question why Fermat wasn't more explicit about the nature of
his~$E$ is an interesting one.  Note that Fermat was involved in an
acrimonious rivalry with Descartes.  Descartes thought that one of the
strengths of \emph{his} own method was that it was purely algebraic.
It is possible that Fermat did not wish to elaborate on the meaning
of~$E$ because he wished to avoid criticism by Descartes.

\section{Methodological issues in 17th century historiography}
\label{method}

On 15 august 1657, Fermat sent a letter to Kenelm Digby (1603--1665).
The letter was entitled ``Remarques sur l'arithm\'etique des infinis
de S. J. Wallis''.  The letter contains a critique of Wallis's
infinitesimal method that reveals as much about Fermat's own position
as about Wallis's method.  As we will see, certain aspects of Wallis's
method \emph{not} criticized by Fermat are as interesting as the
actual criticisms.

\subsection{Fermat's letter to Digby}

The letter is cited by A.~Malet \cite[p.~37, footnote~48]{Mal96}.  It
is also mentioned by J.~Stedall, who goes on to say that
\begin{quote}
The [mathematical] details of the subsequent argument need not concern
us here (Stedall \cite[p.~12]{St01}).
\end{quote}
We will be precisely interested in the mathematical, as well as the
``metamathematical'' issues involved.  Fermat summarizes his objection
to Wallis's method in the following terms:
\begin{quote}
Mais, de m\^eme qu'on ne pourroit pas avoir la raison de tous les
diam\`etres pris ensemble des cercles qui composent le c\^one \`a ceux
du cylindre circonscrit, si on n'avoit la quadrature du triangle; non
plus que la raison des diam\`etres des cercles qui composent le
cono\"\i de parabolique \`a ceux qui font le cylindre circonscrit, si
on n'avoit la quadrature de la parabole ; ainsi on ne pourra pas
conno\^\i tre la raison des diam\`etres de tous les cercles qui
composent la sph\`ere \`a ceux des cercles qui composent le cylindre
circonscrit, si l'on n'a pas la quadrature du cercle (Fermat
\cite[p.~348]{Fe57}).

\end{quote}
Fermat is making a remarkable claim to the effect that in order to
find the quadrature of the circle, Wallis is exploiting the quadrature
of the circle itself.  Fermat appears to be criticizing an alleged
circularity in Wallis's reasoning.  

Apart from the issue of the potency of his critique, what is striking
about it is the aspect of Wallis's method that Fermat is \emph{not}
criticizing.  Namely, what emerges from Fermat's presentation is that
Fermat is taking the infinitesimal technique itself for granted.  In
the paragraph preceding the one cited above, Fermat talks about
spheres and cylinders being composed of \emph{infinite} families of
parallel circles as a routine matter:

\begin{quote}
D'o\`u il conclut que, puisqu'on a trouv\'e aussi la raison de la
sph\`ere au cylindre circonscrit, ou celle de l'infinit\'e des cercles
parall\`eles, dont on peut concevoir que la sph\`ere est compos\'ee,
\`a pareille multitude de ceux qui se peuvent feindre au cylindre, on
pourra aussi esp\'erer de pouvoir d\'ecouvrir la raison des
ordonn\'ees en la sph\`ere ou au cercle \`a celles du cylindre ou
quarr\'e, savoir la raison des diam\`etres des cercles infinis qui
composent la sph\^ere aux diam\`etres des cercles du cylindre. Ce qui
seroit avoir la quadrature du cercle (Fermat \cite[p.~347-348]{Fe57}).
\end{quote}
Thus, it is not the infinitesimal method itself that Fermat is
criticizing, but rather the logic of Wallis's reasoning.

%This text may have been overlooked by \Stromholm{} (1968, \cite{Strom})
%when he confessed his frustration with apparent lack of any discussion
%of infinitesimals in Fermat.

\subsection{Huygens and Rolle}

C.~Huygens \cite{Hu} declared in 1667 at the French Academy of
Sciences that Fermat's ``$e$'' was an infinitely small quantity:
\begin{quote}
Or, en prenant~$e$ infiniment petite, la m\^eme \'equation donnera la
valeur de~$EG$ lorsqu'elle est \'egale \`a~$EF$ \ldots Ensuite on
divise tous les termes par~$e$ et on d\'etruit ceux qui, apr\`es cette
division, contiennent encore cette lettre, puisqu'ils repr\'esentent
des quantit\'es infiniment petites par rapport \`a ceux qui ne
renferment plus~$e$ (Huygens \cite{Hu}, cited in Trompler and No\"el
\cite[p.~110]{TN}).
\end{quote}
Huygens's interpretation is testimony to the enduring influence of
Fermat's method of adequality already in the 17th century.  Yet,
Huygens may have been putting in Fermat's mouth words that did not
emanate therefrom.%
\footnote{See also Section~\ref{con} for Leibniz's view on Fermat's
method.}

Similarly, Michel Rolle in 1703 claimed a connection between Fermat's
$a$ and~$e$ and Leibniz's~$dx$ and~$dy$:
\begin{quote} 
En 1684, Mr de Leibniz donna dans des journaux de Leipzig des exemples
de la formule ordinaire des tangentes, et il imposa le nom
d'\'egalit\'e diff\'erentielle \`a cette formule [\ldots] Mr de
Leibniz n'entreprend point d'expliquer l'origine de ces formules dans
ce projet, ni d'en donner la d\'emonstration [\ldots] Au lieu de l'$a$
\& de l'$e$, il prend~$dx$ \&~$dy$ (Rolle 1703, \cite[p.~6]{Ro}).
\end{quote}
Yet Rolle was an enemy of the calculus, and his identification of~$e$
and~$dx$ may have been due to his eagerness to denigrate Leibniz.  How
are we to avoid this type of pitfall in analyzing Fermat's oeuvre?  In
discussing Fermat's mathematics, two traps are to be avoided:

\begin{enumerate}
\item
Whiggish history, that is, ``the study of the past with direct and
perpetual reference to the present" (H.~Butterfield 1931
\cite[p.~11]{Bu}).
%
% -- in our case, reading Fermat not in his own terms
%but as a precursor of Leibniz and Newton.
%
A convincing reading of Fermat must be solidly grounded in the 17th
century and its ideas, rather than 19th or 20th centuries and their
ideas.
\item
One needs to consider the possibility that Fermat wasn't working with
clear concepts that were destined to become those of the calculus.
\end{enumerate}

We will discuss each of them separately in Subsections~\ref{whig}
and~\ref{clear}.

\subsection{Whig history}
\label{whig}

As far as trap (1) is concerned, it was precisely the risk of
tendentious re-writing of mathematical history that prompted Mancosu
to observe that
\begin{quote}
the literature on infinity is replete with such `Whig' history.
Praise and blame are passed depending on whether or not an author
might have anticipated Cantor and naturally this leads to a completely
\emph{anachronistic} reading of many of the medieval and later
contributions (Mancosu \cite[p.~626]{Ma09}) [emphasis added--the
authors].
\end{quote}
Thus, Cauchy has been often presented anachronistically as a sort of
proto-Weierstrass.  Such a Cauchy--Weierstrass tale has been
critically analyzed by B\l aszczyk et al. \cite{BKS}, Borovik et
al.~\cite{BK}, Br\aa ting \cite{Br}, and Katz \& Katz (\cite{KK11a},
\cite{KK11b}).  Whiggish tendencies in Leibniz scholarship were
analyzed by Katz \& Sherry (\cite{KS1}, \cite{KS2}).

To guard against this trap, we will eschew potential 19th and 20th
century ramifications of Fermat's work, and focus entirely on its 17th
century context.  More specifically, we will examine a possible
connection between Fermat's adequality and Leibniz's Transcendental
Law of Homogeneity (TLH), which play parallel roles in Fermat's and
Leibniz's approaches to the problem of maxima and minima.  Note the
similarity in titles of their seminal texts: \emph{Methodus ad
Disquirendam Maximam et Minimam} (Fermat, see Tannery
\cite[pp.~133]{Tan}) and \emph{Nova methodus pro maximis et minimis
\ldots} (Leibniz \cite{Le84}).

Leibniz developed the TLH in order to mediate between assignable and
inassignable quantities.  The TLH governs equations involving
differentials.  H.~Bos interprets it as follows:
\begin{quote}
A quantity which is infinitely small with respect to another quantity
can be neglected if compared with that quantity.  Thus all terms in an
equation except those of the highest order of infinity, or the lowest
order of infinite smallness, can be discarded.  For instance,
\begin{equation}
\label{adeq2}
a+dx =a
\end{equation}
\[
dx+ddy=dx
\]
etc.  The resulting equations satisfy this [\dots] requirement of
homogeneity (Bos \cite[p.~33]{Bos} paraphrasing Leibniz 1710
\cite[p.~381-382]{Le10b}).
\end{quote}
The title of Leibniz's 1710 text is \emph{Symbolismus memorabilis
calculi algebraici et infinitesimalis in comparatione potentiarum et
differentiarum, et de lege homogeneorum transcendentali}.  The
inclusion of the transcendental law of homogeneity 
%
%(\emph{lege homogeneorum transcendentali} in the ablative case,
%
(\emph{lex homogeneorum transcendentalis})
%
%in the nominative case), 
%
in the title of the text attests to the importance Leibniz attached to
this law.

The ``equality up to an infinitesimal'' implied in TLH was explicitly
discussed by Leibniz in a 1695 response to Nieuwentijt, in the
following terms:

\begin{quote}
Caeterum \emph{aequalia} esse puto, non tantum quorum differencia est
omnino nulla, sed et quorum differentia est incomparabiliter parva; et
licet ea Nihil omnino dici non debeat, non tamen est quantitas
comparabilis cum ipsis, quorum est differencia (Leibniz 1695
\cite[p.~322]{Le95}) [emphasis added--authors]
\end{quote}
Translation:
\begin{quote}
Besides, I consider to be \emph{equal} not only those things whose
difference is entirely nothing, but also those whose difference is
incomparably small: and granted that it [i.e., the difference] should
not be called entirely Nothing, nevertheless it is not a quantity
comparable to those whose difference it is.
\end{quote}

How did Leibniz use the TLH in developing the calculus? The issue can
be illustrated by Leibniz's justification of the last step in the
following calculation:
\begin{equation}
\label{41c}
\begin{aligned}
d(uv) &= (u+du)(v+dv)-uv=udv+vdu+du\,dv \\ & =udv+vdu.
\end{aligned}
\end{equation}

The last step in the calculation~\eqref{41c}, namely
\[
{udv+vdu} + {du\,dv} = {udv+vdu}
\]
is an application of Leibniz's TLH.  In his 1701 text {\em Cum
Prodiisset\/} \cite[p.~46-47]{Le01c}, Leibniz presents an alternative
justification of the product rule (see Bos \cite[p.~58]{Bos}).  Here
he divides by~$dx$ and argues with differential quotients rather than
differentials.  The role played by the TLH in this calculation is
similar to that played by adequality in Fermat's work on maxima and
minima.

\subsection{Clear concepts?}
\label{clear}

Was Fermat working with clear concepts that were destined to become
those of the calculus?  This is a complex question, that in fact
conflates two separate issues: (a) were Fermat's ideas clear? and (b)
were Fermat's ideas destined to become those of the calculus?  Even
the latter formulation is questionable, as the definite article in
front of ``calculus'' disregards that fact, emphasized by
H.~Bos~\cite{Bos}, that the principles of Leibnizian calculus based on
differentials differ from those of modern calculus based on functions.
Thus the answer to (b) is certainly ``we don't know'', though there is
a parallelism between adequality and TLH as we argued in
Subsection~\ref{whig}.  As far as question (a) is concerned, it needs
to be pointed out that such concerns are as old as the critique of
Fermat's method by Descartes, who precisely thought that Fermat was
confused and his method in the category of a lucky guess.%
\footnote{\label{mah2}See e.g., Mahoney \cite[pp.~180--181]{Mah94}.}
However, most modern scholars don't share Descartes' view.  Thus,
H.~Breger wrote:
\begin{quote}
brilliant mathematicians usually are not so very confused when talking
about their own central mathematical ideas [\ldots] I would like to
stress that my hypothesis renders Fermat's mathematics clear and
intelligible, that the hypothesis is supported by several philological
arguments, and that it does not need the assumption that Fermat was
confused (Breger \cite[pp.~193--194]{Bre94}).
\end{quote}
For all our disagreements with Breger, this is one point we can agree
upon.

\section{Did Fermat make a mistake?}
\label{three}

In interpreting Fermat's adequality, one has to keep in mind that
certain crucial components of the conceptual structure of the calculus
were as yet unknown to Fermat.  A striking illustration of this is
what \Stromholm{} refers to as Fermat's ``mistake'' (\Stromholm{}
\cite[p.~51]{Strom}).  Is it really a mistake, a redundancy, or
neither?

In this section, we will write~$A$ and~$E$ in place of~$a$ and~$e$,
following \Stromholm.  After forming the difference~$f(A)-f(A+E)$ and
cancelling out terms not containing~$E$, Fermat writes (in Tannery's
French translation):
\begin{quote}
On divisera tous les termes par~$E$, ou par une puissance de~$E$
\cite[p.~121]{Fer}.%
\footnote{``We divide all the terms by~$E$, or by some power of~$E$''
(the latter clause corresponds to Step~\eqref{higher} of Fermat's
procedure as outlined in Subsection~\ref{summary}).}
\end{quote}
\Stromholm~comments as follows concerning this division:
\begin{quote}
Fermat told his readers that one was to divide by {\em some\/} power
of~$E$.  This, of course, was wrong as can be seen from
\[
f(A+E)-f(A) = \sum_{n=1}^\infty \frac{E^n}{n!} f^{(n)}(A).
\]
Still, his mistake was understandable (\Stromholm{}
\cite[p.~51]{Strom}).
\end{quote}
\Stromholm~does not explain how exactly it can be seen from the Taylor
series expansion that Fermat was wrong.  \Stromholm~continues as
follows:
\begin{quote}
As he could not possibly foresee the peculiarities and future
significance of the~$f^{(n)}(A)$, he guarded himself against the
possibility that~$f'(A)$ be zero, a case which might conceivably (to
him) turn up in some future problem'' (ibid.).
\end{quote}
In the framework of \Stromholm's narrow interpretation of the method,
no such case could conceivably turn up (but see below), and for a
reason unrelated to any ``future significance of~$f^{(n)}(A)$''.
Namely, if the derivative vanishes identically (note that Fermat's
calculation is a symbolic manipulation without assigning a particular
value to the variable~$A$), then the original function itself is
identically constant by the fundamental theorem of calculus (FTC).
The latter wasn't proved until the 1670s, by Isaac Barrow.  Not being
aware of the FTC, Fermat was apparently also unaware of the fact that
no such future problem could turn up, and therefore left in the phrase
``some power of~$E$'', even though only the first power is relevant.

Have we shown then that Fermat's description of his method contains a
mistake, or at least a redundancy?  This is in fact not the case.
Giusti (2009, \cite[\S 6]{Giu}) notes that in the fifth article
\emph{Appendice \`a la M\'ethode du Maximum et Minimum} there is an
example involving radicals.  In this case, the method of adequality as
applied by Fermat leads to an expression in which (before division)
the least power of~$E$ is~$2$ rather than~$1$, and one does have to
divide by~$E^2$ instead of by~$E$.%
\footnote{\label{mah}\Stromholm's mistake is repeated by Mahoney in
both editions of his book: ``in the problems Fermat worked out, the
proviso of repeated division by [$E$] was unnecessary'' (Mahoney
\cite[p.~165]{Mah94}).}

What about the above argument based on the FTC then?  What happened is
that Fermat performs a series of algebraic simplifications so as to
eliminate the radicals from the equation, a point also noted by
Andersen~\cite[p.~59]{An}.  The result in this case is an expression
where the least power of~$E$ is~$2$.  Such manipulations are not
mentioned in the algorithmic/formal description of the method of
adequality, and were not taken into account in \Stromholm's
description of Fermat's ``mistake''.

This example is a striking illustration of the fact that Fermat's
method of adequality is not a single method but rather a cluster of
methods.  The algorithmic procedure described by Fermat at the outset
is merely a kernel common to all applications, but in each application
the kernel is applied somewhat differently.  We will analyze one such
difference in the next two sections.

\section{Comparing the first and second problems in {\em M\'ethode\/}}
\label{comparing}

In this section, we will analyze some apparent dissimilarities between
Fermat's approaches to the first two problems in his {\em M\'ethode
pour la recherche du maximum et du minimum\/} \cite[p.~122--123]{Fer}.

The first problem involves splitting a segment into two subsegments so
as to maximize the area of the rectangle formed by the subsegments as
sides.  The solution is a straightforward application of the formal
algorithmic technique he outlined on the previous page
\cite[p.~121]{Fer}.

Fermat's second problem involves finding the equation of a tangent
line to a parabola.  Mathematically speaking, the second problem is
equivalent to the first.  Namely, given a point, say~$P$, on the
parabola, we would write down the point-slope formula for a line
through~$P$ (with the point fixed and slope, variable), form the
difference between the point-slope and the formula for the parabola,
and look for an extremum of the resulting expression.

But did Fermat view it that way?  He does not appear to have described
it that way.  Forming the difference of the two formulas is an
algebraic procedure.  Did Fermat have such an algebraic procedure in
mind, or would such an approach go beyond the geometric framework as
it actually appears in Fermat?

What Fermat did write is that the point on the tangent line lies {\em
outside\/} the parabola.  Having stated this geometric fact, Fermat
proceeds to write down an {\em inequality\/} expressing it.  Is the
resulting inequality an inessential embellishment of this particular
application of the method of maxima and minima, or is it an essential
part of the argument?

At least on the surface of it, Fermat's formulation is unlike the
earlier case where one obtains an adequality immediately, due to the
nature of the problem, without using an intermediate {\em
inequality\/}.

The passage from the inequality to adequality, depending on whether it
is seen as an essential ingredient in the argument, may or may not
make the second example different from the first one, as we discuss in
Subsection~\ref{four}.

\subsection{Tangent line and convexity of parabola}
\label{four}

Consider Fermat's calculation of the tangent line to the parabola, see
\cite[p.~122-123]{Fer}.  To simplify Fermat's notation, we will work
with the parabola~$y=x^2$, or
\[
\frac{x^2}{y}=1.
\]
To understand what Fermat is doing, it is helpful to
think of the parabola as a level curve of the two-variable
function~$\frac{x^2}{y}$.

Given a point~$(x,y)$ on the parabola, Fermat wishes to find the
tangent line through the point.  Fermat exploits the geometric fact
that by convexity, a point
\[
(p,q)
\]
on the tangent line lies {\em outside\/} the parabola.  He therefore
obtains an inequality equivalent in our notation to~$\frac{p^2}{q}>1$,
or~$p^2>q$.  Here~$q=y-e$, and~$e$ is Fermat's magic symbol we wish to
understand.  Thus, we obtain
\begin{equation}
\label{41}
\frac{p^2}{y-e}>1.
\end{equation}
At this point Fermat proceeds as follows:
\begin{enumerate}
\item[(i)]
\label{shalosh}
he writes down the inequality~$\frac{p^2}{y-e}>1$, or~${p^2}>{y-e}$;
\item[(ii)]
\label{arba}
he invites the reader to {\em ad\'egaler\/} (to ``adequate'');
\item[(iii)]
\label{ve}
he writes down the
adequality~$\frac{x^2}{p^2}=_{\text{AD}}^{\phantom{I}}\frac{y} {y-e}$;
\item[(iv)] 
he uses an identity involving similar triangles to substitute 
\[
\frac{x}{p}=\frac{y+r}{y+r-e}
\]
where~$r$ is the distance from the vertex of the parabola to the point
of intersection of the tangent to the parabola at~$y$ with the axis of
symmetry,
\item[ {(v)}] he cross multiplies and cancels identical terms on right
and left, then divides out by~$e$, discards the remaining terms
containing~$e$, and obtains the solution~$y=r$.%
\footnote{In Fermat's notation~$y=d$,~$y+r=a$. Step (v) can be
understood as requiring the
expression~$\frac{y}{y-e}-\frac{(y+r)^2}{(y+r-e)^2}$ to have a double
root at~$e=0$, leading to the solution~$y=r$ or in Fermat's
notation~$a=2r$.}
\end{enumerate}

What interests us here are the steps~(i) and (ii).  How does Fermat
pass from an inequality to an adequality?  Giusti already noted that
\begin{quote}
Comme d'habitude, Fermat est autant d\'etaill\'e dans les exemples
qu'il est r\'eticent dans les explications.  On ne trouvera donc
presque jamais des justifications de sa r\`egle des tangentes (Giusti
\cite{Giu}).
\end{quote}
In fact, Fermat provides no explicit explanation for this step.
However, he uses the same principle of applying the defining relation
for a curve to points on the tangent line to the curve.  Note that
here the quantity~$e$, as in~$q=y-e$, is positive: Fermat did not have
the facility we do of assigning negative values to variables.  Thus,
\Stromholm~notes that Fermat
\begin{quote}
never considered negative roots, and if~$A=0$ was a solution of an
equation, he did not mention it as it was nearly always geometrically
uninteresting (\Stromholm{} \cite[p.~49]{Strom}).
\end{quote}

Fermat says nothing about considering points~$y+e$ ``on the other
side'', i.e.  further away from the vertex of the parabola, as he does
in the context of applying a related but different method, for
instance in his two letters to Mersenne (see \cite[p.~51]{Strom}), and
in his letter to Br\^ulart~\cite{Fer2}.%
\footnote{This was noted by Giusti \cite{Giu}.}
Now for positive values of~$e$, Fermat's inequality~\eqref{41} would
be satisfied by a {\em transverse ray\/} (i.e., secant ray) starting
at~$(x,y)$ and lying outside the parabola, just as much as it is
satisfied by a tangent ray starting at~$(x,y)$.  Fermat's method
therefore presupposes an additional piece of information, privileging
the tangent ray over transverse rays.

\subsection{Two interpretations of the geometric ingredient}

What is the nature of the additional piece of information that would
privilege the tangent ray?  There are two possible approaches here:
\begin{itemize}
\item
one can argue that the additional piece of information is derived from
the geometric context: namely, the tangent line provides a better
approximation than a transverse line, motivating the passage to an
adequality.
\item
one can argue that both Fermat's geometric context (tangent line being
outside the parabola) and his inequality ($p^2>q$) are merely
incidental, and that Fermat's procedure here is purely algebraic,
namely, equivalent to forming the difference between the formula for a
line through~$(x,y)$ and the formula for the parabola, and seeking an
extremum as before.
\end{itemize}

In support of the geometric interpretation, it could be said that, in
order to understand Fermat's procedure, other than treating it as a
magician's rabbit out of a hat, we need to relate to the geometric
context.  To passage from inequality to adequality therefore only
becomes intelligible as something less mysterious, if one assumes
that~$e$ is {\em small\/} and exploits the geometric background with a
better rate of approximation provided by the tangent line as compared
to a transverse ray.  Str\o mholm \cite{Strom} similarly emphasizes
the role of the {\em smallness\/} of~$e$ in Fermat's thinking.  To
assert that Fermat's procedure using the symbol~$e$ is purely
formal/algebraic in the context of this particular example, is to
assert that Fermat is a magician, not a mathematician.

In support of the algebraic interpretation, it could be said that
Fermat writes ``ad\'egalons donc, d'apr\`es la m\'ethode
pr\'ec\'edente'' with reference to the second example, apparently
implying that the method of the first and second examples is
comparable.  Since both methods contain a common kernel as we
discussed in Section~\ref{three}, the reference to the previous
example is not conclusive.

Treating the geometric ingredient as an essential part of the proof in
this case is the more appealing option, suggesting that the method of
tangents is not a direct application of the kernel of the method of
maxima and minima, but rather exploits additional geometric
information in a crucial way.  Similarly, K.~Pedersen comments:
\begin{quote}
The inequality~$IO>IP$ holds for all curves concave with respect to
the axis, and the inequality~$IO<IP$ for convex curves.  For curves
without points of inflection it is possible from these inequalities to
find a magnitude depending on~$a-e$ and~$x-e$ which has an extreme
value for~$x-e=x$ (Pedersen%
\footnote{This is the same author as Andersen \cite{An}.}
1980 \cite[p.~28]{Pe80}).
\end{quote}
She continues:
\begin{quote}
Neither in his {\em Methodus\/} nor in Fermat's later writings,
however, is there any indication that this was the way he related his
method of tangents to his method of maxima and minima (ibid.),
\end{quote}
and concludes that
\begin{quote}
Descartes was right after all in raising the objection that the method
of tangents was not a direct application of the method of maxima and
minima (ibid.).
\end{quote}

We saw that the geometric content of the argument dealing with
tangents to parabolas tends to go counter to the formal/algebraic
interpretation.  Barner \cite[p.~34]{Bar} attempts to save the day by
declaring that Fermat made a mistake.  According to Barner, the point
we denoted~$(p,q)$ in Section~\ref{four} should not be outside the
parabola at all, but rather should be on the parabola.  The line
passing through~$(p,q)$ should not be the tangent line, but rather a
transverse line, whereas~$e$ is the difference in the ordinates of the
two points on the parabola.  In this way, the {\em inequality\/}
should not be there in the first place, but should rather be an {\em
equality\/}.

In fact, what Barner describes with uncanny accuracy is the method of
Beaugrand, probably dating from 1638, as analyzed by Str\o mholm
\cite[p.~64-65]{Strom}.  We therefore find Barner's explanation, as
applied to Fermat, to be forced (see also Subsection~\ref{barner}).
The formal/algebraic interpretation of adequality is unconvincing in
this particular case.

\section{Fermat's treatment of the cycloid}
\label{cycloid}

The cycloid is generated by marking a point on a circle and tracing
the path of the point as the circle rolls along a horizontal straight
line.  If the marked point is the initial point of contact of the
circle with the line, and the circle rolls to the right, then the
ordinate%
\footnote{In this case the \emph{ordinate} refers to the horizontal
coordinate.}
of the marked point is given by the difference of the length of arc
traversed (the distance the center of the circle has moved) and the
distance of the point from the vertical line through the center of the
circle.%
\footnote{Assuming the circle to have radius~$1$, the equation of the
cycloid as described above is~$x=\theta - \sin \theta,\quad y=1 - \cos
\theta$.}

\subsection{Fermat's description}

\begin{figure}%[ht]
\includegraphics[height=3in]{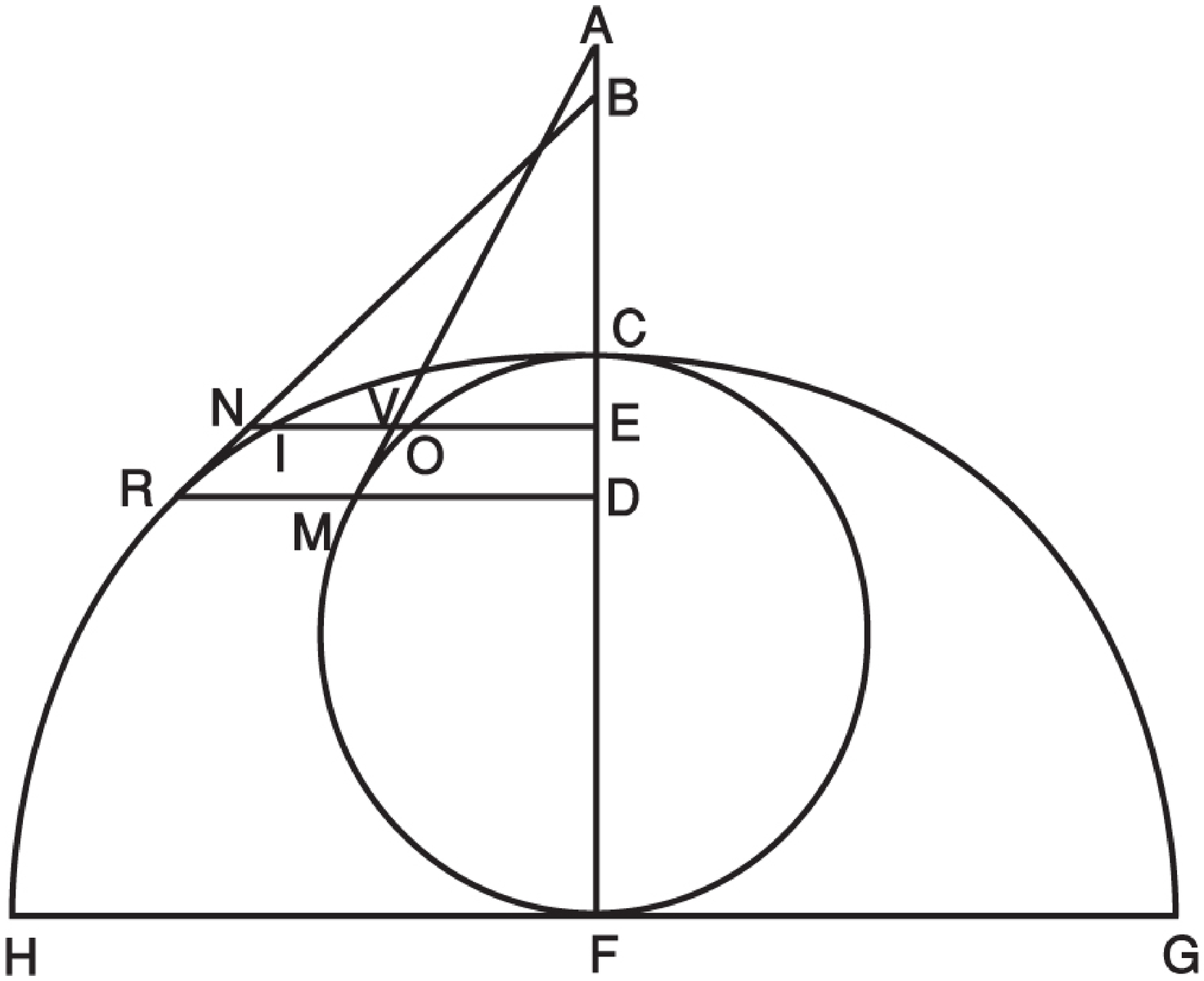}
\caption{\textsf{Fermat's cycloid}}
\label{apr12}
\end{figure}

Fermat's description of the cycloid is based on a diagram \cite[Figure
103, p. 163]{Tan} reproduced in Figure~\ref{apr12}.  Let~$R$ be a
point on the cycloid and~$D$ the point of intersection of the
horizontal line~$\ell$ through~$R$ with the axis of symmetry of the
cycloid generated by one full revolution of the circle. If~$M$ is the
point of intersection of~$\ell$ with the generating circle when
centered on the axis of symmetry, and~$C$ is the apex of that circle
then in the words of Fermat:
\begin{quote}
La propri\'et\'e sp\'ecifique de la courbe est que la droite
$RD$ est \'egale \`a la somme de l'arc de cercle~$CM$ et de
l'ordonn\'ee~$DM$%
\footnote{To compare Fermat's description with the parametric
description given in the previous footnote, we note that length of the
segment~$RD$ is~$\pi-x=\pi-\theta + \sin\theta$, while~$\pi-\theta$ is
the length of the arc~$CM$, and the length~$DM$ equals~$\sin\theta.$}
(\cite[p~144]{Fer}).
\end{quote}

Let~$r$ be the tangent line to the cycloid at~$R$, and~$m$ the tangent
line to the circle at~$M$.  To determine the defining relation of~$r$,
Fermat considers the horizontal line~$NIVOE$ passing through a
point~$N\in r$.  Here~$I$ is the first point of intersection with the
cycloid, while~$V$ is the point of intersection with~$m$, and~$O$ is
the point of intersection with the generating circle, and~$E$ is the
point of intersection with the axis of symmetry.

The defining relation for~$r$ is derived from the defining relation
for the cycloid using adequality.  The defining relation for the
point~$I$ on the cycloid is
\[
IE=OE + {\rm arc}\,CO.
\]
By \emph{adequality}, Fermat first replaces~$I$ by the point~$N\in r$:
\begin{equation}\label{adeq}
NE =_{\rm AD} OE +{\rm arc}\,CO=OE+{\rm arc}\, CM - {\rm arc}\, MO.
\end{equation}
Then Fermat replaces~$O$ by the point~$V\in m$, and the~$\text{arc}\,
MO$, by the length of the segment~$MV\subset m$.  This produces the
linear relation
\[
NE =_{\rm AD}VE+{\rm arc}\,CM - MV,
\]
yielding the equation of the tangent line~$r$ to the cycloid at~$R$ as
a graph over the axis of symmetry.  The distance~$NE$ is expressed in
terms of the distance~$VE$, where~$V\in m$, and the distance~$MV$
along that tangent line.  Thinking in terms of slope relative to the
variable distance~$DE$ (which corresponds to the parameter~$e$ in the
example of the parabola), Fermat's equation says that the slope of~$r$
relative to~$DE$ is the slope of~$m$ minus the proportionality factor
of~$MV$ relative to~$DE$.%
\footnote{The slope of the tangent line relative to the axis of
symmetry, or equivalently, relative to the~$y$ axis, given by
elementary calculus is~$ \frac{d(\pi-x)}{d\theta} /
\frac{dy}{d\theta}= \frac{-1}{\sin\theta} +\frac
{\cos\theta}{\sin\theta}.$ The length~$MV$ equals~$e/(\sin\theta)$ and
the slope of the tangent line to the circle relative to the~$y$ axis
is~$\frac{\cos\theta}{\sin\theta}$, in agreement with Fermat's
equation.}
To summarize, Fermat exploited two adequations in his calculation:
\begin{enumerate}
\item the length of a segment along~$m$ adequals the length of a
segment of a circular arc, and
\item 
the distance from the axis of symmetry to a point on~$r$ (or~$m$)
adequals the distance from the axis to a corresponding point on the
cycloid (or circle).
\end{enumerate}
As Fermat explains,
\begin{quote}
Il faut donc {\em ad\'egaler} (\`a cause de la propr\'et\'e
sp\'ecifique de la courbe qui est \`a consid\'erer sur la tangente)
cette droite~$\frac{za-ze}a$ [i.e.,~$NE$] \`a la somme~$OE + {\rm
arc}\, CO$ \dots [et] d'apr\`es la remarque pr\'ec\'edente,
substituer, \`a~$OE$, l'or\-donn\'ee~$EV$ de la tangente, et \`a
l'arc~$MO$, la portion de la tangente~$MV$ qui lui est adjacente
(Fermat \cite[p.~144]{Fer}; \cite[p.~228]{Tan}).
\end{quote}
The procedure can be summarized in modern terms by the following
principle.

\begin{principle}
The tangent line to the curve is defined by using adequality to
linearize the defining relation of the curve, or ``\emph{ad\'egaler}
(\`a cause de la propr\'et\'e sp\'ecifique de la courbe qui est \`a
consid\'erer sur la tangente)."
\end{principle}

Fermat uses the same argument in his calculation of the tangents to
other transcendental curves whose defining property is similar to the
cycloids and involves arc length along a generating curve.  For a
discussion of some of these examples, see Giusti \cite{Giu} and
Itard~\cite{It}.

\subsection{Breger's interpretation}
\label{52}

Breger claims that Fermat
\begin{quote}
does not make use of an approximate equality between arc length of the
circle and a segment of the tangent. The latter interpretation usually
is made, but it plainly contradicts the text: Fermat explicitly calls
the straight line~$DE$ ``recta utcumque assumpta" (Fermat 1891, 163)
[droite arbitraire], that is,~$DE$ is not infinitely small or ``very
small" \cite[p.~206]{Bre94}.
\end{quote}
Is it convincing to argue, as Breger does, that the hypothesis of
smallness plays no role here?  It is certainly true that the axis of
symmetry contains points arbitrarily far from the fixed point~$D$.
However, being able meaningfully to apply the defining relation for
the curve to the tangent line, is contingent upon the fact that for
points near the point of tangency the tangent line gives a good
(second order) approximation to the curve.

Why does Fermat's calculation give the tangent to the cycloid?  In
modern terms, this is because of the quadratic order in~$e$ for the
error term which results from making each of the following two
substitutions: (1)~substituting the point~$N$ on the tangent line for
the point~$I$ on the cycloid and (2) substituting the length of a
segment of the tangent to the circle for arc~$CO$ on the circle.

%The vanishing of the expression at~$e=0$, mentioned by Breger, is
%nsufficient, and second order approximation in~$e$ is indispensable.

Breger's analysis of the example proceeds as follows: 
\begin{quote}
Fermat (1891, 163 last paragraph) looks at the
expression~$NE-OE-CM+MO~$ ($CM$ and~$MO$ being arcs of the
circle)[eq. (\ref{adeq})].  This expression takes a minimum, namely
zero, if the point~$E$ coincides with the point~$D$. Then Fermat
replaces~$MO-OE$ by~$MU-UE$ [$U=V$ above]. Now if the point~$E$ is not
too far away of [sic] the point~$D$, then
$$ MO-OE\geq -MD$$ as well as
$$ MU -UE\geq -MD$$ and equality holds if and only if the points~$E$
and~$D$ coincide. Therefore~$NE -UE -CM+MU$ has a minimum, namely
zero, if the points~$E$ and~$D$ coincide (Breger
\cite[p.~206]{Bre94}).%
\footnote{The segment~$MD$ in Breger's text does not appear in either
expression~$NE- OE - CM + MO$ or~$NE - UE - CM + MU$ and is clearly a
misprint for~$MC$.}
\end{quote}

Note that Breger suppresses the term~$e=DE$ by setting it equal to
zero (in discussing the expression~$NE -UE -CM+MU$) at
stage~\ref{cancel}, namely prior to division by~$e$,%
\footnote{\label{meaningless2}See footnote~\ref{meaningless1}.}
and with this concludes his argument in favor of a ``minimum''
interpretation.  Such a procedure is however meaningless and certainly
cannot be attributed to Fermat.  Fermat clearly states:
\begin{quote}
Divisons par~$e$; comme il ne reste ici aucun terme superflu, il n'y a
pas d'autre suppression \`a faire (Fermat \cite[p.~145]{Fer}).
\end{quote}
Fermat thus asserts that there are no terms to be suppressed after
division by~$e$ in the example of the cycloid.  Breger's erroneous
suppression of~$DE$ prior to division by~$e$ is not accidental, but
rather stems from a desire to force a ``minimum'' interpretation on
Fermat.

Breger concludes his discussion of the cycloid by pointing out that
\begin{quote}
It is hard for the modern reader to get rid of the \emph{limit ideas}
in his mind, and so he considers the replacement of a very small arc
length of the circle by a very small segment of the tangent to be
quite natural.  But this is \emph{not the way of Fermat} (Breger 1994
\cite[p.~206]{Bre94}) [emphasis added--the authors].
\end{quote}
Granted, modern ``limit ideas'' as applied to Fermat are a tell-tale
sign of Whiggish history discussed in Section~\ref{method}.  However,
setting~$e$ equal to zero prior to division by~$e$ is hardly ``the way
of Fermat'', either.

To conclude, Fermat's treatment of the cycloid and other
transcendental curves cannot be accounted for by means of a
formal/algebraic reading, and requires an element of approximation for
a convincing interpretation.

\subsection{Who erred: Fermat or Barner?}
\label{barner}

In a recent article, K.~Barner (2011, \cite{Bar}) expresses agreement
with Breger, while noting an additional fine point.  He claims that,
while \emph{aequare} and \emph{adaequare} are semantically equivalent,
Fermat uses the term \emph{adaequare} when defining a dependence
relation between two variables.  Thus, the adequation
$f(x+h)=_{AD}f(x)$ defines the dependence of~$x$ on~$h$.  Barner goes
on to interpret Fermat's method in terms of (Dini's) implicit function
theorem.  

On page~23, Barner writes: ``Und genau dies ist die Bedeutung des
Wortes \emph{ad\ae quare}: es bezeichnet die Gleichheit zwischen zwei
Termen, die \emph{keine Identit\"at} sondern eine von der Identit\"at
verschiedene \emph{Relation} definiert'' [And this is exactly the
meaning of the word \emph{adaequare}: it indicates the equality
between two terms, an equality which defines not an identity, but
rather a relationship distinct from identity].

Barner further claims that there is not really one uniform method, and
to assume that Fermat has a clearly formulated method is an
anachronism: ``Zu unterstellen, daÃ dies f\"ur Fermat ebenfalls eine
selbstverst\"andliche Routine gewesen sei, \emph{das} ist ein
\emph{Anachronismus}'' [To assume that for Fermat this [method] was a
self-evident/obvious procedure is an anachronism].

On page 27, Barner goes on to claim that Fermat had to use a bit of
`trickery' because his method is not completely correct: ``Ich
m\"ochte einfach nur verstehen, warum Fermat mit seiner Trickserei
Erfolg hat, obwohl der Ansatz, den er dabei macht, die Idee seiner
`Methode' nicht ganz korrekt wiedergibt'' [I would simply like to
understand why Fermat was successful with his little trickery, even
though the assumption that he makes does not represent the basic idea
of his `method' with full accuracy].

Barner's approach ultimately leads him astray when dealing with
tangents.  Barner seeks to replace the tangent line by a secant line
(page 34), since his interpretation forces him to assume that the
second point near the point of tangency is on the parabola itself
(point~$P$ in his diagram on page 34) rather than on the tangent line
(point~$O$ in his diagram), so as to get an exact equality: ``Aber der
Punkt~$P$ (und damit bei Fermat auch~$O$!) liegt damit auf der Kurve
und erlaubt es ihm, \emph{suivant la propri\'et\'e sp\'ecifique de la
ligne courbe}, die algebraischen Eigenschaften der vorgelegten Kurve
auch f\"ur den Punkt~$O$ in Anspruch zu nehmen'' [But the
point~$P$--and thus, according to Fermat, the point~$O$ as well--lies
thereby on the curve and allows him, \emph{suivant la propri\'et\'e
sp\'ecifique de la ligne courbe}, to claim the algebraic properties of
the curve under consideration for the point~$O$ as well].

How plausible is Barner's claim?  Barner's claim contradicts Fermat's
explicit statement, which Barner had just quoted:
\begin{quote}
Apr\`es avoir donn\'e le nom, tant \`a notre parall\`ele qu'\`a tous
les autres termes de la question, tout le m\^eme qu'en la parabole, je
consid\`ere derechef cette parall\`ele, comme si le point qu'elle a
dans la tangente \'etoit en effet en la ligne courbe, et suivant la
propriet\'e sp\'ecifique de la ligne courbe, je compare cette
parall\`ele par ad\'egalit\'e avec l'autre parall\`ele tir\'ee du
point donn\'e \`a l'axe ou diam\`etre de la ligne courbe.  Cette
comparaison par ad\'egalit\'e produit deux termes in\'egaux qui enfin
produisent l'\'egalit\'e (selon ma m\`ethode), qui nous donne la
solution de la question (Fermat quoted in Barner
\cite[pp.~35-36]{Bar}).
\end{quote}

Thus, Barner has pursued his interpretation to the point of
contradicting Fermat's own comments.  Fermat writes explicitly that he
is applying the defining property of the curve to points on the
\emph{tangent} line.  He says it again in his discussion of the
cycloid as we have shown (see Subsection~\ref{52}).

On page~36, Barner goes on to say that he can't understand why Fermat
kept making the same mistake: ``Warum hat Fermat sein
widerspr\"uchliches Vorgehen bei allen seinen Beispielen zur
Tangentenmethode immer wieder erneut verwendet? Warum hat er die
Sekante, mit der er de facto operiert, nie erw\"ahnt? Ich wei\ss{} es
nicht'' [Why did Fermat continually repeat his inconsistent procedure
for all his examples for the method of tangents?  Why did he never
mention the secant, with which he in fact operated? I do not know].

Barner's dilemma does not arise in our interpretation, which avoids
attributing errors to Fermat.

\section{Bachet's semantic calque}
\label{calque}

The choice of the Latin verb \emph{adaequo} in Bachet's translation of
Diophantus is explained in Bachet's notes as follows:
\begin{quote}
Since in questions of this kind, Diophantus nearly%
\footnote{\emph{proxime} in the Latin.}
equ\-ates the sides of the squares%
\footnote{Mathematically speaking, it would apparently have been more
correct to write ``squares of the sides''.}
that are being sought, to some side, but he does not properly equate
them, he calls this comparison {\parisotes\/} and not
'\!\!$\iota\sigma\acute{o}\tau\eta\varsigma$.  We too call it not
equality but adequality, just as we also translate
$\pi\acute\alpha\rho\iota\sigma o\nu$ as adequal.%
\footnote{``Quia enim in huiusmodi quaestionibus Diophantus, cuidam
lateri adaequat proxime latera quadratorum quaesitorum, non autem
aequat proprie, vocat ille hanc comparitionem
$\pi\alpha\rho\iota\sigma\acute{o}\tau\eta\tau\alpha$ non autem \,
'\!\!$\iota\sigma\acute{o}\tau\eta\tau\alpha$.  Nos etiam non
aequalitatem sed adaequalitatem appellamus, sicut etiam
$\pi\acute{\alpha}\rho\iota\sigma o \nu$ vertimus adaequale.''}
\end{quote}

Here Bachet clearly differentiates the meanings of \emph{aequo} and
\emph{adaequo}: the former is ``equal'', the latter ``nearly equal''.%
\footnote{See a related discussion in Itard
(1974,~\cite[p.~338]{It74}).}
Notice that what Bachet performs here is a semantic calque: the Greek
\emph{para} and \emph{iso\=o} are individually translated to \emph{ad}
and \emph{aequo}, and recombined to produce \emph{adaequo}.  Bachet
does not have the Latin meaning of \emph{adaequo} in mind, but rather
a new meaning derived from his understanding of the term coined by
Diophantus.

\section{Breger on mysteries of {\em adaequare\/}}

In 1994, Breger sought to challenge what he called the ``common dogma"
to the effect that

\begin{quote}
{\em Fermat uses ``adaequare'' in the sense of ``to be approximately
equal" or ``to be pseudo-equal" or ``to be counterfactually equal"\/}
\cite[p.~194]{Bre94} [emphasis in the original--the authors].
\end{quote}

After some introductory remarks related to the dating and editing of
various manuscripts, he continues
\begin{quote}
Having made these introductory remarks I want to put forward my {\em
hypothesis: Fermat used the word ``adaequare" in the sense of ``to put
equal\/}'' \cite[p.~197]{Bre94} [italics in the original].
\end{quote}

In this section, we compare Breger's interpretation and that of the
viewpoint of \Stromholm~and Giusti that there are (at least) two
different approaches in Fermat.  One is based on the insight from
Pappus, and involves a symmetric relation between two equal values
near an extremum.  The other is based on the insight from Diophantus
and the method of \parisotes, which exploits an approximation.  In the
first approach, ``adequality" has an operational meaning of ``setting
equal" as Breger contends.  However, in the second approach (see
Section~\ref{cycloid}), as well as in the application of the method to
tangents (see Section~\ref{comparing}), the element of approximation
is essential.
 
%Revise: Delete the following subsection. Diophantus's method
%is described in the section devoted to the mathematics of
%parisotes, problems 12,14, and 17.
%\subsection{The mathematical background}

%It may be helpful briefly to summarize the mathematical background of
%Diophantus' problem in modern terms.  The underlying principle is that
%if there is a rational point on a variety then in fact rational points
%are dense.  In other words, there are rational points arbitrarily
%close to a given irrational one~$(x,y)$.  Diophantus' problem is
%equivalent to the case of the specific irrational point
%\[
%x=y=\sqrt{\tfrac{2a+1}{2}}.
%\]
%Since irrational points are outside of Diophantus' conceptual
%framework, Diophantus can only formulate the problem in terms of
%seeking rational~$x$ and~$y$ that are sufficiently close to each other
%(in the sense of the bound~\eqref{42} below), and satisfy
%\begin{equation}
%x^2+y^2=2a+1.
%\end{equation}

\subsection{The philology of \parisotes} 
\label{philo}

Breger presents several arguments in defense of his hypothesis. The
first two arguments are based on dictionary definitions, first of the
Latin {\em adaequare\/} and second of the original Greek term
{\parisotes}. Both of these arguments are flawed.

Breger questions why Fermat would choose to employ the term {\em
adaequare\/} in a sense different from that given in standard Latin
dictionaries:
\begin{quote}
There are well established Latin words to indicate an approximate
equality, namely ``approximare" or, more frequently used in classical
Latin as well as in the 17$^{\text{th}}$ century, ``appropinquare".
It is well known that Fermat's knowledge of Latin and Greek was
excellent, and so if he wanted to tell us that there was an
approximate equality, why should he use a word indicating an equality?
As far as I know, this question has not been answered nor even
discussed by any adherent of the dogma of Fermat interpretation
(Breger \cite[p.~198]{Bre94}).
\end{quote}
The first point is that there was no pre-existing word meaning ``to
set two things approximately equal".  Non-mathematical texts do not
talk that way; when Polybius says, for example, that Rome and Carthage
were~$\pi\acute\alpha\rho\iota\sigma\alpha$ in their power
(Polyb. 1.13.8), nobody takes him to be stating a mathematical
equivalence. And the words that Breger suggests, \emph{approximare} or
\emph{appropinquare}, do not mean what Fermat wants to do: he does not
want to ``approximate"~$f(a)$, but rather to \emph{compare} two
expressions,~$f(a)$ and~$f(a+e)$.  If what he wanted to do was to
treat them as being approximately equal, neither the Romans nor the
mathematicians of Fermat's time had a term available for that; and
since the Latin prefix \emph{ad-} commonly translated the Greek prefix
\hbox{\emph{par(a)-},} the term \emph{adaequalitas} would be the
obvious equivalent for \parisotes.  

To give an example of a Greek \emph{par} transformed into a Latin
\emph{ad}, note that the rhetorical figure that the Greeks called
\emph{paronomasia} ($\pi\alpha\rho o\nu o
\mu\alpha\sigma\acute\iota\alpha$), whereby a speaker makes a pun on
two similar but not identical words, was called in Latin
\emph{adnominatio}.  But as far as Fermat is concerned, the obvious
reason that he used \emph{adaequalitas} is that that was the term that
Bachet used to translate the Greek \parisotes.  Morever, Fermat
himself explicitly refers to the Greek term.

The fact that Breger is aware of this problem is evident when he
writes:
\begin{quote}
Therefore Xylander and Bachet were right in using the word
``adaequalitas" in their translations of Diophantus [\ldots] although
that does not imply that they had understood the mathematics of the
passage (Breger \cite[p.~200]{Bre94}).
\end{quote}
To argue his hypothesis, Breger is led to postulate that Xylander and
Bachet misunderstood Diophantus!  Breger appears to acknowledge
implicitly that Bachet's intention was to engineer a semantic calque
(see Section~\ref{calque}), but argues that Bachet's calque was a vain
exercise, to the extent that \emph{adaequo} already has the meaning of
\emph{pariso\=o}, albeit not the meaning Bachet had in mind.

\subsection{Greek dictionary}

Breger's ``argument [that] is just based on the Greek dictionary"
\cite[p. 199]{Bre94} is misconceived.  The Parisian \emph{Thesaurus
Graecae Linguae} of 1831-75 \cite{TGL2}%
\footnote{Breger gives the dates 1842-7, understating the magnitude of
the task.}  
on which Breger bases himself, was a reissue of Henricus Stephanus'
dictionary of the same name \cite{TGL1}.  The latter, a work of
stupendous scholarship, was published in 1572, three years before
Xylander's original edition of Diophantus.  Neither Stephanus'
Thesaurus, nor the cheaper and hence much more widely available
pirated abridgement of it by Johannes Scapula \cite{TGL3} (originally
printed in 1580 and reprinted innumerable times afterward), nor
Stephanus' 1582 reissue, included the word \parisotes, though both
Stephanus and Scapula did have~$\pi\acute\alpha\rho\iota\sigma
o\varsigma$, \emph{aequalis}, \emph{uel compar}: ``equal or similar"
(with the note that the writers on rhetoric used the expression to
mean \emph{prope aequatum}, ``made nearly equal") and other words
formed from~$\pi\alpha\rho\acute\alpha$ and '\!\!$\acute\iota\sigma
o\varsigma$.

In the nineteenth-century Paris reissue, the term \parisotes{} was
added by Karl Wilhelm Dinsdorf, one of the editors, who cited
Diophantus and translated ``\AE qualitas."
%
%\footnote{ Dinsdorff's translation of \parisotes{} as \emph{equality}
%is at tension with our textual analysis of Diophantus in
%Section~\ref{eight}, which suggests that \emph{approximation} is
%inherent in its meaning.  The meaning of \emph{equality} may have been
%the result of a back-translation from the traditional meaning of the
%Latin \emph{adaequo}, along the lines of Breger's approach.}
%
The Paris Stephanus, however, is not today ``the best Greek
dictionary", contrary to Breger's claim.  The most (and in fact, for
the time being, the only) authoritative dictionary today is
Liddell-Scott-Jones \cite{LSJ}, which defines \parisotes{} as
``approximation to a limit, Dioph.~5.17."  But in fact all Greek and
Latin dictionaries (with the exception of the original Stephanus,
which was chiefly the result of Stephanus' own scholarship and that of
Guillaume Bud\'e) are secondary sources, recording the meanings that
others have given to the words; so both Dindorf and
Liddell-Scott-Jones were simply recording the translations then
current for Diophantus' use of \parisotes.  Both ``equality" and
``near equality" are possible meanings for \parisotes; and it is up to
the editors and critics of Diophantus to tell the lexicographers which
meaning he intended, not the other way around.

Fermat is very unlikely to have used either the hugely expensive and
rare Stephanus or the easily available Scapula, when Bachet had
printed Diophantus with a facing Latin translation; and his care in
saying that he is using the term \emph{adaequentur ut loquitur
Diophantus}, ``as Diophantus says", seems to make it clear that he is
specifically \emph{not} using the Latin term in its usual meaning of
``be set equal", but, rather, in a meaning peculiar to Diophantus.
That meaning, as we argue in Section~\ref{eight}, can only be ``be set
approximately equal", as it was correctly understood by Bachet.

\section{The mathematics of \parisotes}
\label{eight}

In this section we will analyze the problems in which Diophantus
introduces the term {\parisotes}.  

\subsection{The \parisotes{} of Diophantus}
\label{61}

The relevant problems are in Book Five, problems 12, 14, and 17 in
Bachet's Latin translation.%
\footnote{The numbering is the same in Wertheim's German translation.
However, in Ver Eecke's French translation, the corresponding problems
are 9, 11, and~14, which is the numbering referred to by Weil.}
Problem~14 uses the term \parisos~(see \cite{Bachet}).

The term \parisos{} occurs on page 310 of Bachet's Latin translation
\emph{Diophanti Alexandrini, Arithmeticorum Liber V}.  Here the left
column is the Greek original, while the right column is the Latin
translation.  The term \parisos{} is the first word on line 8 in the
left column.  In the right column, line 13, we find its Latin
translation \emph{adaequalem}.

\subsection{Lines 5-8 in Diophantus}
In more detail, lines 5-8 in the Greek original contain the following
phrase containing \parisos:

%\vfill\eject %\section{hello}

%
%dei oun ton i. dielein eis treis tetragonous opos ekastou tetragonon e
%pleura parisos e m ia^stigma
%

\newcommand{\stigma}{\text{\includegraphics[height=6pt]{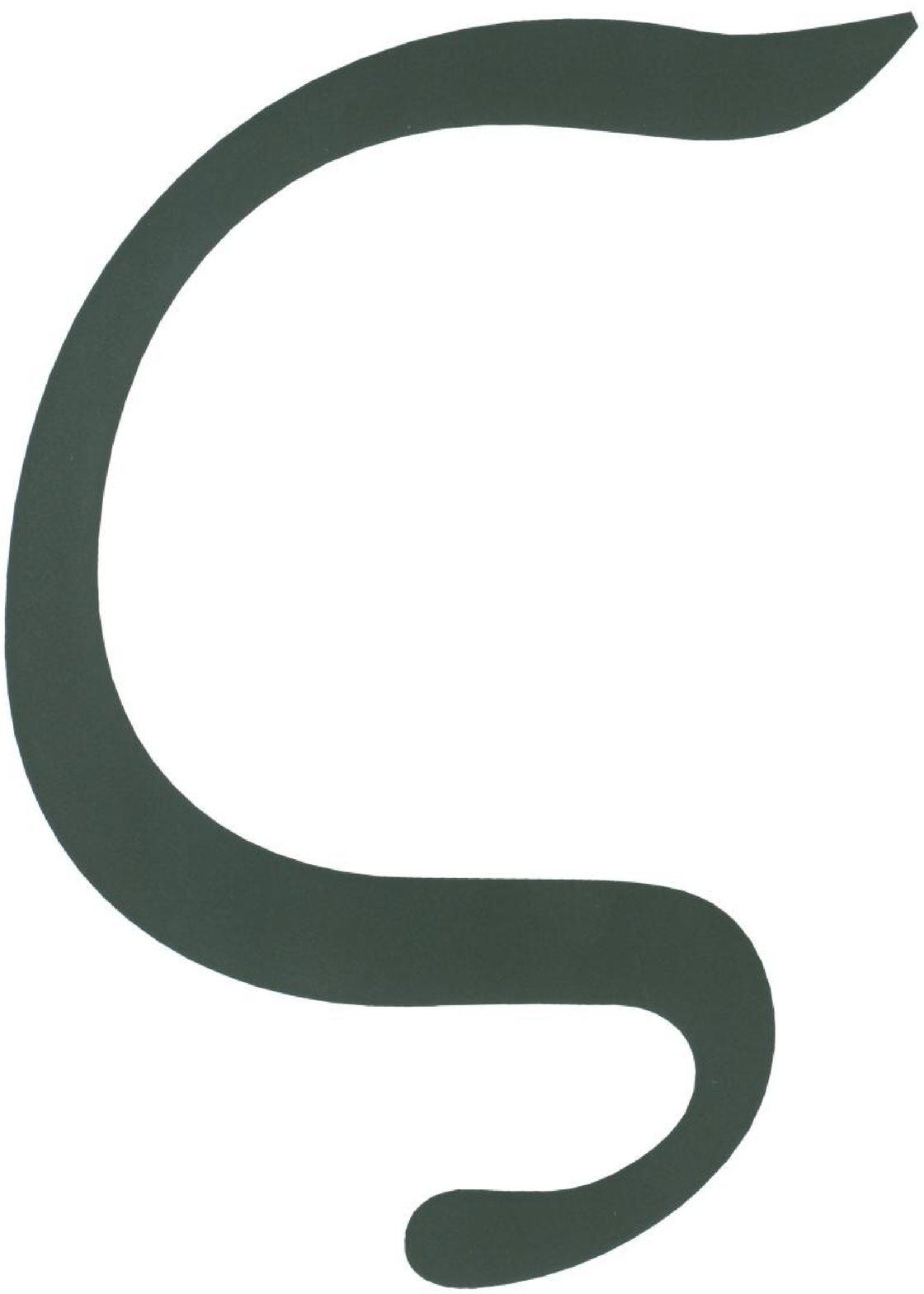}}}

\begin{quote}
$\delta\varepsilon\tilde\iota$
%dei
$o \tilde{\stackrel{,}{\upsilon}} \nu$
%
%oun
%
$\tau \grave o\nu\; \bar \iota.$
%
%ton i
%
$\delta\iota\varepsilon\lambda\varepsilon\tilde\iota\nu$
%
%dielein
%
$\epsilon \! \stackrel{,}{\iota} \! \varsigma$
%
%eis
%
$\tau\rho\varepsilon\tilde\iota\varsigma$
%
%treis
%
$\tau\varepsilon\tau\rho\alpha\gamma\acute\omega\nu
o\upsilon\varsigma$
%
%tetragonous
%
`\!$\acute o\pi\omega\varsigma$
%
%opos
%
\\
%
%ekastou
%
`\!$\varepsilon\kappa\acute\alpha\sigma\tau o\upsilon$
%
%tetragonou
%
$\tau\varepsilon\tau\rho\alpha\gamma\acute\omega\nu o\upsilon$
%
%e
%
`\!$\eta$
%
%pleura
%
$\pi\lambda\varepsilon\upsilon\rho\grave\alpha$
%
%parisos
%
$\pi\acute\alpha\rho\iota\sigma{o}\varsigma$
$\underset{`}{\tilde\eta}$
$\bar\mu'$
$\bar\iota\alpha^{\bar
%
%
%\varsigma
\stigma
}.$%
\footnote{\label{stigma}The letter appearing in the exponent under the
bar is the letter \emph{stigma} (a ligature of~$\sigma$ and~$\tau$).}
\end{quote}
\newcommand{\longs}{\text{\hskip1.3pt\includegraphics[height=9pt]{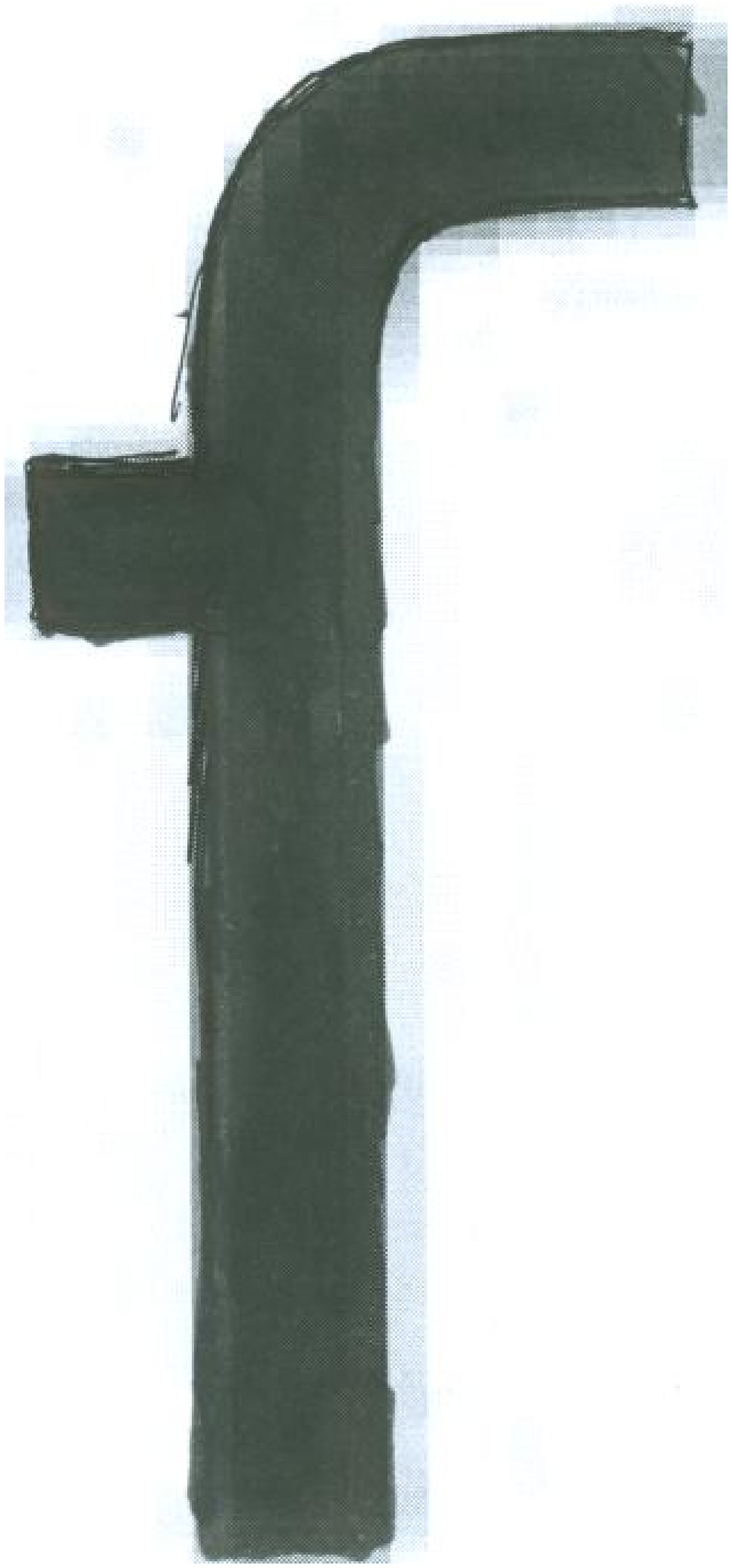}}}
Latin translation as printed:%
\footnote{Here all three occurrences of the symbol \longs{} are long
``s''s.}
\begin{quote}
Oportet igitur diuidere 10. in tres quadratos, vt vniu\longs
cuiu\longs que quadrati latus \longs it ad\ae quale
vnitatibus~$\frac{11}{6}$.
\end{quote}
Latin translation in modern letters:
\begin{quote}
Oportet igitur dividere 10 in tres quadratos, ut uniuscuiusque
quadrati latus sit adaequale unitatibus~$\frac{11}{6}$.
\end{quote}
English translation:
\begin{quote}
So we have to divide ten into three squares, so that the side of each
square is adequal to~$\frac{11}{6}$ units.
\end{quote}

\subsection{Lines 14-16 in Diophantus}

Lines 14-16 in the Greek original contain the following phrase
containing \parisos:

\begin{quote}
$\delta\varepsilon\tilde\iota$
%dei
$o \tilde{\stackrel{,}{\upsilon}} \nu$
%
%oun
%
$\tau \grave o\nu$
%
%ekaste
%
`\!$\varepsilon\kappa\acute\alpha\sigma\tau \eta$
%
%ton
%
$\tau\tilde\omega\nu$
%
%pleuron
%
$\pi\lambda\varepsilon\upsilon\rho\tilde\omega\nu$
%
%touton
%
$\tau o\acute\upsilon\tau\omega\nu$
%
%paraskeuasai
%
$\pi\alpha\rho\alpha\sigma\kappa\varepsilon\upsilon
\acute\alpha\sigma\alpha\iota$~$\pi\acute\alpha\rho\iota\sigma o\nu$
$\bar\iota\alpha^{\bar\stigma}.$%
\footnote{See footnote~\ref{stigma}.}
\end{quote}
Latin translation as printed:
\begin{quote}
Oportet igitur horum cuiu\longs uis lateri ad\ae qualem
facere~$\frac{11}{6}$.
\end{quote}
Latin translation in modern letters:
\begin{quote}
Oportet igitur horum cuiusvis lateri adaequalem facere~$\frac{11}{6}$.
\end{quote}
English translation:
\begin{quote}
So the side of each of these we have to make adequal
to~$\frac{11}{6}$.
\end{quote}

\subsection{Lines 20-21 in Diophantus}
Lines 20-21 in the original Greek:
\begin{quote}
$\delta\epsilon\tilde\iota$
%
%oun
% 
$o \tilde{\stackrel{,}{\upsilon}} \nu$
%
%ekasten
%
`\!$\varepsilon\kappa\acute\alpha\sigma\tau\eta\nu$ 
%
%pleuran
%
$\pi\lambda\varepsilon\upsilon\rho\grave\alpha\nu$
$\kappa\alpha\tau\alpha\sigma\kappa\varepsilon\upsilon
\acute\alpha\sigma\alpha\iota$~$\overline{\nu}\varepsilon$.
\end{quote}
Latin as printed:
\begin{quote}
Oportet itaque vnumquodque latus ad\ae quare ip\longs i 55.
\end{quote}
Latin in modern letters:
\begin{quote}
Oportet itaque unumquodque latus adaequare ipsi 55.
\end{quote}
English:
\begin{quote}
So each side we have to make adequal to 55.
\end{quote}

\subsection{Line 25 in Diophantus}
Line 25: Greek:
\begin{quote}
%
%tauta
%
$\tau\alpha\tilde\upsilon\tau\alpha$
%
%isa
%
'\!$\acute\iota\sigma\alpha$
%
%monasi
%
$\mu o\nu\acute\alpha\sigma\iota$~$\overline{\iota}$.
\end{quote}
Latin as printed:
\begin{quote}
h\ae c \ae quantur 10.
\end{quote}
Latin in modern letters:
\begin{quote}
haec aequantur 10.
\end{quote}
English:
\begin{quote}
These equal 10 units.
\end{quote}

\subsection{\parisos{} as approximate equality}

Diophantus wishes to represent the number 10 as a sum of three
squares.  The solution he eventually finds is denoted~$(\alpha, \beta,
\gamma)$ by Ver Eecke \cite[p.~205]{Ver}.  In this notation,
Diophantus seeks~$\alpha, \beta, \gamma$ each of which is as close as
possible (\parisos) to the fraction~$\frac{11}{6}$.  The solution
eventually found
is~$\alpha=\frac{1321}{711}$,~$\beta=\frac{1288}{711}$,
$\gamma=\frac{1285}{711}$.

Thus, Diophantus specifically uses the term \parisos{} to describe the
way in which the numbers~$\frac{1321}{711}$,~$\frac{1288}{711}$,
and~$\frac{1285}{711}$ approximate the fraction~$\frac{11}{6}$.
Interpreting Diophantus' \parisotes{} as anything other than
``approximate equality" is therefore purely whimsical.  In the next
subsection we describe the method of Diophantus in more detail.

\subsection{The method of Diophantus}

The following description of the method which Diophantus called
{\parisotes\/} is based on the notes of Ver Eecke~\cite{Ver},
Wertheim~\cite{Wer}, and the paper of Bachmakova \cite{Bach}.
Diophantus does not have a theory of equations, but gives an algorithm
for solving a class of problems by solving a particular example.  We
will explain his method in terms of equations for the reader's
convenience.  In all three examples involving {\parisotes}, the
problem is to express a certain number as the sum of two or three
rational squares, with a further inequality constraining the
individual values.

In problem 12, one seeks two rational squares whose sum is 13 and each
one is greater than~6.

In problem 14, one seeks three rational squares whose sum is 10 and
such that each one is greater than~3.  

In problem 17, one seeks two rational squares whose sum is 17 and such
that each one is less than 10.

In all three cases, the first step is to find a rational square
approximately equal to either one half or one third of the desired
sum, depending on whether two or three numbers are sought.  In
addition, Diophantus uses the existence of a preliminary partial
solution involving two or three rational squares giving the desired
sum,~$N$, but failing to satisfy the required inequalities. Call these
numbers~$a_i$ and denote the rational square close to one half or one
third of the desired sum by~$b$. He then expresses~$b$ as a sum
$b=a_i+ c_i$, and considers an equation of the form
\[
\Sigma (a_i +t c_i)^2 =N=\Sigma a_i^2.
\]
Canceling equal terms on the two sides and dividing by~$t$ gives a
linear equation for~$t$ with a rational solution~$t_0$.  Since~$t=1$
is nearly a solution, the exact solution~$t_0$ is close to~$1$ and
defines the numbers~$a_i+ t_0 c_i$ satisfying the required inequality.

As mentioned in Subsection~\ref{61}, Diophantus' \parisos{} refers to
the approximate equality of each of the~$a_i+t_0 c_i$ to the original
fraction~$b$.  Ver Eecke explains the matter as follows in a
footnote:

\begin{quote}
$\pi\alpha\rho\iota\sigma\acute{o}\tau\eta\tau {o}\varsigma$%
\footnote{\label{f19}The spelling of \parisotes{} changed because of
declension.  It appears here in the genitive case.}
%
%this footnote is deleted in the published text at David's request
%
%parisotos
%
$\alpha\!\!^{^{,}}\gamma \omega \gamma \tilde \eta$, c'est-\`{a}-dire
``la voie de la quasi-\'{e}galit\'{e}" ou d'approximation vers une
limite; m\'ethode dont l'objet est de r\'{e}soudre des probl\`{e}mes
tels que celui de trouver 2 ou 3 nombres carr\'{e}s, dont la somme est
un nombre donn\'{e}, et qui sont respectivement au plus rapproch\'{e}s
d'un m\^{e}me nombre (Ver Eecke \cite[p.~203, footnote 2]{Ver}).
\end{quote}
The method, as explained here, is remarkably similar to Fermat's.  It
starts with a near equality and defines a quadratic equation which is
solved as an exact equality.  The method of {\parisotes} might be
called ``the method of nearby values" (see Subsection~\ref{ety}).
G.~Lachaud \cite[Section~5]{La} similarly speaks of \parisotes{} as
involving approximation.

\subsection{Breger on Diophantus}

Breger acknowledges that
\begin{quote} 
At the first step - and this is the relevant step characterizing the
method of \parisotes{} - he finds a positive rational {\em z} with the
property~$2z^2\approx 2a+1$ \cite[p.~200]{Bre94}.
\end{quote}
After discussing the Diophantine problems, Breger claims that
\begin{quote}
the approximate equality which in fact occurs in Diophantus's
Arithmetic is only due to the fact that Diophantus seeks rational
solutions \cite[p.~202]{Bre94}.
\end{quote}
Breger's claim that Diophantus is seeking an approximate rational
solution is in error.  In fact, Diophantus starts from an approximate
rational solution, to derive an equation for an \emph{exact} rational
solution.  The issue of approximation is not related to the
distinction rational/irrational since irrational solutions were
outside the scope of Diophantus's conceptual framework.  Breger
further remarks that
\begin{quote} 
[i]t is strongly misleading to mix Fermat's notation with our own and
to describe his method in these cases by something like ``$f(A)$
adaequatur~$f(A-E)$'', as is often done \cite[p.~204]{Bre94}.
\end{quote}
How is one to interpret Breger's comment?  The fact is that setting
the former expression adequal to the latter is precisely what Fermat
does in the very first example, using~$a+e$ rather than~$a-e$.  In the
French translation which uses modern notation (in place of Fermat's
original notation \`a la Vi\`ete for the mathematical expressions),
the phrase appears in the following form:
\begin{quote} 
Soit maintenant~$a+e$ le premier segment de~$b$, le second
sera~$b-a-e$, et le produit des segments:
\newline~$ba-a^2 +be-2ae-e^2$;
\newline Il doit \^{e}tre ad\'{e}gal\'{e} au pr\'{e}c\'{e}dent:
$ba-a^2$ (Fermat \cite[p.~122]{Fer}).
\end{quote}
In the original Latin version, the last two lines above read as
follows:
\[
B\,{\rm in }\,A -Aq. +B\,{\rm in }\,E - A\,{\rm in }\,E\,{\rm
bis}-Eq.,
\]
\hfil\text{quod debet adaequari superiori rectangulo}
\[
B\,{\rm in}\,A -Aq. 
\]
(see Fermat \cite[p.~134]{Tan}).  Breger continues, ``The method
consists in finding a polynomial in~$E$ which takes a minimum
if~$E=0$" \cite[p.~204]{Bre94}.  However, Fermat does not merely apply
his method to polynomials.  He ultimately applied the method in far
greater generality:
\begin{quote}
At first Fermat applies the method only to polynomials, in which case
it is of course purely algebraic; later he extends it to increasingly
general problems, including the cycloid (Weil \cite[p~.1146]{We}).
\end{quote}
In Section~\ref{cycloid} we argued that transcendental curves such as
the cycloid necessarily require an element of approximation.

\subsection{The Diophantus--Fermat connection}
The mathematical areas Diophantus and Fermat were working in were
completely different.  It was arithmetic and number theory in the case
of Diophantus, and geometry and calculus in the case of Fermat.  Such
a situation creates a fundamental problem: if one rejects the
``approximation'' thread connecting Diophantus to Fermat, why exactly
did Fermat bring Diophantus and his terminology into the picture when
working on problems of maxima and minima?

Breger's solution to the problem is to declare that Diophantus was
talking about {\em minima\/} and Fermat was also talking about
minima.  What kind of minima was Diophantus talking about?  Breger's
answer is that Diophantus' minimum is\ldots~$0$.  Namely,~$|x^2-y^2|$
is always bigger than zero, but gets arbitrarily close to it:
\begin{quote}
The minimum evidently is achieved by putting~$x$ equal to~$y$, and
that is why the method received its name ``method of \parisotes{} or
putting equal".  As there is a minimum idea in the Diophantus passage,
Fermat's reference to Diophantus becomes intelligible (Breger
\cite[p.~201]{Bre94}).
\end{quote}
Does Fermat's reference to Diophantus become intelligible by means of
the observation that~$|x^2-y^2|$ gets arbitrarily close to~$0$?
According to Breger's hypothesis, Diophantus was apparently led to
introduce a new term, \parisotes, to convey the fact that every
positive number is greater than zero.  To elaborate on Breger's
hypothesis, zero is the \emph{infimum} of all positive numbers, which
is presumably close enough to the idea of a \emph{minimum}.  All this
is supposed to explain the connection to Fermat's method of minima.

There are at least two problems with such a reading of Diophantus.
First, did Diophantus have the number zero?  Second, the condition
that Diophantus imposes is merely the bound
\begin{equation}
\label{42}
|x^2-y^2|<1,
\end{equation}
rather than any stronger condition requiring the expression
$|x^2-y^2|$ to be arbitrarily close to zero.  Thus, Diophantus was
\emph{not} concerned with the infimum of~$|x^2-y^2|$.  If one drops
the approximation issue following Breger, the entire
Diophantus--Fermat connection collapses.

\section{Refraction, adequality, and Snell's law}
\label{nine}

In addition to purely mathematical applications, Fermat applied his
adequality in the context of the study of refraction of light, so as
to obtain Snell's law.  Thus, in his {\em Analyse pour les
r\'efractions\/}, Fermat sets up the formulas for the length of two
segments
\[
\text{CO}.m=\sqrt{m^2n^2+m^2e^2-2m^2be}, \quad
\text{IO}.b=\sqrt{b^2n^2+b^2e^2+2b^2ae},
\]
representing the two parts of the trajectory of light across a
boundary between two regions of different density, and then writes:
\begin{quote}
La somme de ces deux radicaux doit \^etre \emph{ad\'egal\'ee},
d'apr\`es les r\`egles de l'art, a la somme~$mn+bn$ (Fermat
\cite[p.~150]{Fer}) [emphasis added--the authors]
\end{quote}

Fermat explains the physical underpinnings of this application of
adequality in his \emph{Synth\`ese pour les r\'efractions\/}, where he
writes that light travels slower in a denser medium (Fermat
\cite[p.~151]{Fer}).  Fermat states the physical principle
underpinning his mathematical technique in the following terms:
\begin{quote}
Notre d\'emonstration s'appuie sur ce seul postulat que la nature
op\`ere par les moyens et les voies les plus faciles et les plus
ais\'ees  \cite[p.~152]{Fer},
\end{quote}
and goes on to emphasize that this is contrary to the traditional
assumption that ``la nature op\`ere toujours par les lignes les plus
courtes''.  Rather, the path chosen is the one traversed ``dans le
temps le plus court'' (ibid.).  This is Fermat's principle of least
time in optics.  Its implicit use by Fermat in his {\em Analyse pour
les r\'efractions\/} in conjunction with adequality, is significant.

Namely, this physical application of adequality goes against the grain
of the formal/algebraic approach.  The latter focuses on the higher
multiplicity of the root of the polynomial~$f(a+e)-f(a)$ at an
extremum~$a$, where the extremum can be determined by an algebraic
procedure without ever assigning any specific value to~$e$, and
obviating the need to speak of the nature of~$e$.  In this method,
denoted M2 by \Stromholm, the symbol~$e$ could be a formal variable,
neither small or large, in fact without any relation to a specific
number system.

However, when we apply a mathematical method in physics, as in the
case of the refraction principle provided by Snell's law, mathematical
idealisations of physical magnitudes are necessarily numbers.  The
principle is that the light chooses a trajectory~$\tau_0$ which
minimizes travel time from point A to point B.  To study the principle
mathematically is to commit oneself to comparing such a trajectory to
other trajectories~$\tau_s$ in a family parametrized by a numerical
parameter~$s$.  Here we need not assume an \emph{identity} of the line
in physical space with a number line (the hypothesis of such an
identification is called Cantor's axiom in the literature); rather, we
merely point out that a number line is invariably what is used in
mathematical idealisations of physical processes.

In studying such a physical phenomenon, even before discussing the
size of~$e$ (small, infinitesimal, or otherwise), one necessarily
commits oneself to a number system rather than treating~$e$ as a
formal variable.  Fermat's application of adequality to derive Snell's
law provides evidence against a strict formal/algebraic interpretation
of adequality.

Applications to physical problems necessarily involve a mathematical
implementation based on numbers.  Classical physics was done with
numbers, not algebraic manipulations.  One can't model phenomena in
classical physics by means of formal variables.

Certainly when a physicist performs mathematical operations, he does
exactly the same thing as a mathematician does.  However, {\em
modeling\/} physical phenomena by using mathematical idealisation, is
a stage that {\em precedes\/} mathematical manipulation itself.  In
such modeling, physical phenomena get numerical counterparts, and
therefore necessarily refer to a number system rather than formal
variables.

Pedersen \& Pedersen interpret Fermat's deduction of the sine law of
refraction
\begin{quote}
as an early example of the calculus of variations rather than as an
ordinary application of Fermat's method of maxima and minima (Pedersen
\cite{Pe71}),
\end{quote} 
and speculate that Fermat thus anticipated Jacob Bernoulli who is
generally credited with inventing the method of the calculus of
variations in 1696.  Arguably, at least some of the manifestations of
the method of adequality amount to variational techniques exploiting a
small or infinitesimal variation~$e$.

\section{Conclusion}
\label{con}

Should Fermat's~$e$ be interpreted as a formal variable, or should it
be interpreted as a member of an Archimedean continuum ``tending" to
zero?  Or perhaps should adequality be interpreted in terms of a
Bernoullian continuum,%
\footnote{G.~Schubring attributes the first systematic use of
infinitesimals as a foundational concept, to Johann Bernoulli, see
\cite[p.~170, 173, 187]{Sch}.  To note the fact of such systematic use
by Bernoulli is not to say that Bernoulli's foundation is adequate, or
that that it could distinguish between manipulations with
infinitesimals that produce only true results and those manipulations
that can yield false results.  One such infinitesimal distinction was
provided by Cauchy \cite{Ca53} in 1853, thereby resolving an ambiguity
inherent in his 1821 ``sum theorem'' (see Katz \& Katz \cite{KK11b}
for details).}
with~$e$ infinitesimal?  Note that the term ``infinitesimal" was not
introduced until around 1670,%
\footnote{Some contemporary scholars hold that Leibniz coined the term
{\em infinitesimal\/} in 1673 (see Probst \cite{Pr08} and \cite{Mo}).
Meanwhile, Leibniz himself attributed the term to Nicolaus Mercator
(see Leibniz \cite{Le99}).}
so Fermat could not have used it.  Yet infinitely small quantities
were in routine use at the time, by scholars like John Wallis who was
in close contact with Fermat.

While discussions of ``process" are rare in Fermat when he deals with
his~$e$, we mentioned an instance of such use in Subsection~\ref{13}.%
\footnote{See quotation following footnote~\ref{kinetic}.}
Writes \Stromholm:
\begin{quote}
It will not do here to drag forth the time-honoured ``limiting
process" of historians of mathematics [\dots] Fermat was still
thinking in terms of equations; I agree that he stood on the verge of
a period where mathematicians came to accept that sort of process, but
he himself was in this particular case rather the last of the ancients
than the first of the moderns (\Stromholm{} \cite[p.~67]{Strom}).
\end{quote}

In the absence of infinitesimals, there is no possibility of
interpreting {\em smallness\/} other than by means of a process of
{\em tending\/} to zero.  But, as \Stromholm~confirms, such a
discussion is uncharacteristic of Fermat, even at the application
stage, when he applies his method in concrete instances.  Therefore
the infinitesimal interpretation (3) is more plausible than the
kinetic interpretation (2) (see Section~\ref{one}).

To return to the question posed in Section~\ref{one}, as to which of
the three approaches is closest to Fermat's thinking, it could be that
the answer to the riddle is\ldots it depends.  When Fermat presents
his definitional characterisation of adequality, as on the first page
of his {\em M\'ethode pour la recherche du maximum et du minimum\/},
his algorithmic presentation has a strong formal/algebraic flavor.
However, at the application stage, both in geometry and physics, ideas
of approximation or smallness become indispensable.

Breger \cite[205-206]{Bre94} claims that adequality cannot be
interpreted as approximate equality.  Breger's argument is based on
his contention that the Latin term {\em adaequare\/} was not used in
the sense of {\em approximate equality\/} by Fermat's contemporaries.
However, the source of adequality is in the
Greek~$\pi\alpha\rho\iota\sigma\acute{o}\tau\eta\varsigma$
(parisot\=es), rather than the Latin {\em adaequare\/}, undermining
Breger's argument.  The question that should be asked is not whether
Fermat's contemporaries used the term {\em adaequare\/}, but rather
whether they used the infinitely small.  The latter were certainly in
routine use at the time, by some of the greatest of Fermat's
contemporaries such as Kepler and Wallis.

In addition to the 3-way division: formal, kinetic, and infinitesimal,
there is a distinction between (A) Fermat's definition, i.e., synopsis
of the method as it appears in Fermat \cite[p.~121]{Fer}; and (B) what
he actually does when he applies his method.

Fermat's definition (A) does have the air of a kind of a formal
algebraic manipulation.  A formal interpretation of adequality is
certainly {\em mathematically\/} coherent, regardless of what Fermat
meant by it, since one can define differentiation even over a {\em
finite field\/}.  The fact itself of being able to give a purely
algebraic account of this mathematical technique is not surprising.
What is dubious is the claim that at the application stage (B), he is
similarly applying an algebraic procedure, rather than thinking of~$e$
geometrically as vanishing, tending to zero, infinitesimal, etc.

In light of the positivity of Fermat's~$e$ in the calculation of the
tangent line, the formal story would have difficulty accounting for
the passage from inequality to adequality, since the inequality is
satisfied for transverse rays as well as the tangent ray.  To make
sense of what is going on at stage (B), we have to appeal to geometry,
to {\em negligible\/}, {\em vanishing\/}, or {\em infinitesimal\/}
quantities, or their {\em rate\/} or {\em order\/}.

Breger's insistence on the formal interpretation (1), when applied to
the application stage (B), is therefore not convincing.  Fermat may
have presented a polished-up algebraic presentation of his method at
stage (A) that not even Descartes can find holes in, but he gave it
away at stage~(B).

Kleiner and Movshovitz-Hadar note that 
\begin{quote}
Fermat's method was severely criticized by some of his contemporaries.
They objected to his introduction and subsequent suppression of the
mysterious~$e$.  Dividing by~$e$ meant regarding it as not zero.
Discarding~$e$ implied treating it as zero.  This is inadmissible,
\emph{they rightly claimed}.  In a somewhat different context, but
\emph{with equal justification}, Bishop Berkeley in the 18th century
would refer to such~$e's$ as `the ghosts of departed quantities'
\cite[p.~970]{KM} [emphasis added--authors].
\end{quote}
Kleiner and Movshovitz-Hadar feel that Fermat's suppression of~$e$
implies treating~$e$ as zero, and that the criticisms by his
contemporaries and by Berkeley were justified.  However, P.~Str\o
mholm already pointed out in 1968 that in Fermat's main method of
adequality (M1),
\begin{quote}
there was never [\ldots] any question of the variation~$E$ being put
equal to zero.  The words Fermat used to express the process of
suppressing terms containing~$E$ was {\em ``elido''\/}, {\em
``deleo''\/}, and {\em ``expungo''\/}, and in French {\em
``i'efface''\/} and {\em ``i'\^ote''\/}.  We can hardly believe that a
sane man wishing to express his meaning and searching for words, would
constantly hit upon such tortuous ways of imparting the simple fact
that the terms vanished because~$E$ was zero (\Stromholm{}
\cite[p.~51]{Strom}).
\end{quote}

Fermat did not have the notion of the derivative.  Yet, by insisting
that~$e$ is being discarded rather than set equal to zero, he planted
the seeds of the solution of the paradox of the infinitesimal quotient
and its disappearing~$dx$, a century before George Berkeley ever
lifted up his pen to write {\em The Analyst\/}.%
\footnote{\label{discard2}The heuristic procedure of discarding the
remaining terms was codified by Leibniz by means of his Transcendental
Law of Homogeneity (see Section~\ref{whig}).  Centuries later, it was
implemented mathematically in terms of the standard part function,
which associates to each finite hyperreal number, the unique real
number infinitely close to~it.  In 1961, Robinson \cite{Ro61}
constructed an infinitesimal-enriched continuum, suitable for use in
calculus, analysis, and elsewhere, based on earlier work by E.~Hewitt
\cite{Hew}, J.~\Los{} \cite{Lo}, and others.  Applications of
infinitesimal-enriched continua range from aid in teaching calculus
\cite{El, KK1, KK2} to the Boltzmann equation (see
L.~Arkeryd~\cite{Ar81, Ar05}) and mathematical physics (see Albeverio
et al.~\cite{Alb}).  Edward Nelson \cite{Ne} in 1977 proposed an
alternative to ZFC which is a richer (more stratified) axiomatisation
for set theory, called Internal Set Theory (IST), more congenial to
infinitesimals than ZFC.  The hyperreals can be constructed out of
integers (see Borovik et al.~\cite{BJK}).}

After summarizing Nieuwentijt's position on infinitesimals, Leibniz
wrote in 1695:
\begin{quote}
It follows that since in the equations for investigating tangents,
maxima and minima (which the esteemed author [i.e., Nieuwentijt]
attributes to Barrow, although if I am not mistaken Fermat used them
first) there remain infinitely small quantities, their squares or
rectangles are eliminated (Leibniz 1695 \cite[p.~321]{Le95}).
\end{quote}
Leibniz held that methods of investigating tangents, minima, and
maxima involve infinitesimals.  Furthermore, he disagreed with
Nieuwentijt as to the priority of developing these methods,
specifically attributing them to Fermat.  Thus, Leibniz appears to
have felt that Fermat's methods of investigating tangents, minima and
maxima did rely on infinitesimals.  In the absence of explicit
commentary by Fermat concerning the nature of~$E$, Leibniz's view may
be the best 17th century expert view on the matter.

\section*{Acknowledgments} 

We are grateful to Enrico Giusti, Alexander Jones, and David Sherry
for helpful comments.  The influence of Hilton Kramer (1928-2012) is
obvious.

%\newpage

\medskip\noindent \textbf{Mikhail G. Katz} (B.A. Harvard University,
'80; Ph.D. Columbia University, '84) is Professor of Mathematics at
Bar Ilan University.  Among his publications are the following: (with
P.~B\l aszczyk and D.~Sherry) Ten misconceptions from the history of
analysis and their debunking, \emph{Foundations of Science} (online
first); (with A.~Borovik) Who gave you the Cauchy--Weierstrass tale?
The dual history of rigorous calculus, \emph{Foundations of Science}
\textbf{17} ('12), no.~3, 245--276; (with A.~Borovik and R.~Jin) An
integer construction of infinitesimals: Toward a theory of Eudoxus
hyperreals, \emph{Notre Dame Journal of Formal Logic} \textbf{53}
('12), no.~4, 557-570; (with K. Katz) Cauchy's continuum,
\emph{Perspectives on Science} \textbf{19} ('11), no.~4, 426-452;
(with K.~Katz) {Meaning in classical mathematics: is it at odds with
Intuitionism?}  \emph{Intellectica} \textbf{56} ('11), no.~2, 223-302;
(with K.~Katz) {A Burgessian critique of nominalistic tendencies in
contemporary mathematics and its historiography}, \emph{Foundations of
Science} \textbf{17} ('12), no.~1, 51--89; (with K.~Katz) {Stevin
numbers and reality}, \emph{Foundations of Science} \textbf{17} ('12),
no.~2, 109--123; (with E.~Leichtnam) {Commuting and non-commuting
infinitesimals}, to appear in \emph{American Mathematical Monthly};
(with D.~Sherry) {Leibniz's infinitesimals: Their fictionality, their
modern implementations, and their foes from Berkeley to Russell and
beyond}, \emph{Erkenntnis} (online first); (with D.~Sherry) Leibniz's
laws of continuity and homogeneity, \emph{Notices of the American
Mathematical Society}, to appear.

\medskip\noindent \textbf{David M. Schaps} (B.A. Swarthmore, Greek,
highest honors, '67; M.A. Harvard University, Classical Philology,
'70; Ph.D, Classical Philology, Harvard, '72) is Associate Professor
of Classical Studies at Bar Ilan University. Among his publications
are \emph{Economic Rights of Women in Ancient Greece} (Edinburgh '79),
\emph{The Beauty of Japhet} (in Hebrew, Ramat Gan, '89), \emph{The
Invention of Coinage and the Monetization of Ancient Greece} (Ann
Arbor '04), \emph{Handbook for Classical Research} (London/New York
'11), and dozens of articles on ancient Greek history, language, and
literature.

\medskip\noindent \textbf{Steve Shnider} (B.S. Georgetown University
1965, Magna Cum Laude, M.S. Harvard University 1967, Ph.D. Harvard
University 1972) is Professor of Mathematics at Bar Ilan University,
Israel.  He has published over 60 articles and three research
monographs on a wide range of research topics including
C.R. manifolds, gauge theory, supersymmetry, deformation theory,
quantum groups, operads, and most recently Babylonian mathematics.
Among his publications are: Plimpton 322: a Review and a Different
Perspective, \emph{Archive for History of Exact Sciences} \textbf{65}
(2011), 519-566, with John P. Britton and Christine Proust; and the
research monographs: Super Twistor Spaces and Super Yang-Mills
Equation, \emph{Seminaire de Mathematiques Superieures} \textbf{106}
(1991) with Ronny Wells, Jr.; \emph{Quantum Groups, from Coalgebras to
Drinfeld Algebras}, International Press, (1994) with Shlomo Sternberg;
\emph{Operads and their Applications in Algebra, Topology, and
Physics}, Mathematical Surveys and Monographs 96, American
Mathematical Society (2002), with Martin Markl and Jim Stasheff.

\end{document}